\documentclass[10pt]{amsart}
\usepackage[T1]{fontenc}
\usepackage{geometry}
\usepackage[latin1] {inputenc}
\usepackage{amsmath}
\usepackage{amsfonts, amssymb, textcomp}
\usepackage[colorlinks=flase, linkcolor=red,urlcolor=green, citecolor=blue]{hyperref}
\usepackage{subeqnarray}

\usepackage{latexsym}
\usepackage{fancyhdr}
\usepackage{longtable}
\usepackage{amsmath, amssymb}
\usepackage{graphicx}
\usepackage{enumerate}

\numberwithin{equation}{section}

\theoremstyle{plain}
\newtheorem{proposition}{Proposition}[section]
\newtheorem{theorem}[proposition]{Theorem}
\newtheorem{lemma}[proposition]{Lemma}
\newtheorem{corollary}[proposition]{Corollary}
\newtheorem{definition}[proposition]{Definition}

\newtheorem{remark}[proposition]{Remark}

\newcommand{\RR}{\mathbb{R}}
\newcommand{\CC}{\mathbb{C}}
\newcommand{\NN}{\mathbb{N}}

\newcommand{\id}{\operatorname{id}}

\let\on=\operatorname

\newenvironment{demo}[1]{\par\smallskip\noindent{\bf #1.}}{\par\smallskip}

\makeatletter
\@namedef{subjclassname@2020}{%
  \textup{2020} Mathematics Subject Classification}
\makeatother

\newsavebox{\fmbox}
\newenvironment{fmpage}[1]
 {\begin{lrbox}{\fmbox}\begin{minipage}{#1}}
 {\end{minipage}\end{lrbox}\fbox{\usebox{\fmbox}}}

\title[On strong growth conditions for weighted spaces of entire functions]
{On strong growth conditions for weighted spaces of entire functions}

\author[G.~Schindl]{Gerhard Schindl}

\address{G.~Schindl: Fakult\"at f\"ur Mathematik, Universit\"at Wien, Oskar-Morgenstern-Platz~1, A-1090 Wien, Austria.}
\email{gerhard.schindl@univie.ac.at}

\begin{document}

\begin{abstract}
We characterize the inclusion relations between weighted classes of entire functions with rapid decreasing growth and study strong growth comparison relations between given weights. In our considerations first we focus on weights defined in terms of the so-called associated weight function where the weight(system) is based on a given sequence. Then the abstract weight function case is reduced to the weight sequence setting by using the so-called associated weight sequence. Finally, we compare weighted entire function spaces defined in terms of so-called dilatation-type and exponential-type weight systems.
\end{abstract}

\thanks{This research was funded in whole by the Austrian Science Fund (FWF) project 10.55776/P33417}
\date{\today}

\maketitle


\section{Introduction}\label{Introduction}
This article is the direct continuation of the author's recent work \cite{weightedentirecharacterization} and for the notation and ideas used below we refer also to \cite{weightedentirecharacterization}. Weighted spaces of entire functions are defined as follows:
\begin{equation}\label{weightedclassdef}
H^{\infty}_v(\CC):=\{f\in H(\CC): \|f\|_v:=\sup_{z\in\CC}|f(z)|v(|z|)<+\infty\},
\end{equation}
with $H(\CC)$ denoting the class of entire functions and $v:[0,+\infty)\rightarrow(0,+\infty)$ is the (radial and strictly positive) {\itshape weight function.} From now on, except stated explicitly otherwise, we assume that $v$ is
\begin{itemize}
	\item[$(*)$] continuous,
	
	\item[$(*)$] non-increasing and
	
	\item[$(*)$] rapidly decreasing, i.e. $\lim_{t\rightarrow+\infty}t^kv(t)=0$ for all $k\ge 0$.
\end{itemize}
This is the same notion to be a weight (function) as it has been considered in \cite{weightedentirecharacterization} and several other recent papers; see \cite{BonetTaskinen18}, \cite{BonetLuskyTaskinen19}, \cite{solidassociatedweight}, and the informative survey article \cite{Bonet2022survey}. It is known that $H^{\infty}_v(\CC)$ is a Banach space w.r.t. the norm $\|\cdot\|_v$ and similarly such weighted classes can be defined when $\CC$ is replaced by any other (non-empty) open subset $G\subseteq\CC$ and $v: G\rightarrow(0,+\infty)$. Another very prominent example in the literature is the unit disc $G=\mathbb{D}$; in this article we will focus on $G=\CC$ except in Section \ref{equsectforo}.\vspace{6pt}

In \cite{weightedentirecharacterization} we have been concerned with the following problem: Characterize the inclusion relations of weighted entire function classes in terms of the defining weights. More precisely, let $v$ and $w$ be two weights and write $v\hypertarget{ompreceq}{\preceq}w$ if
\begin{equation}\label{bigOrelation}
w(t)=O(v(t))\;\text{as}\;t\rightarrow+\infty.	
\end{equation}
Let us call $v$ and $w$ {\itshape equivalent,} written $v\hypertarget{sim}{\sim}w$, if
$$v\hyperlink{ompreceq}{\preceq}w\;\text{and}\;w\hyperlink{ompreceq}{\preceq}v.$$
Note that relation \hyperlink{sim}{$\sim$} is denoted by the symbol $\approx$ in \cite{Bonet2022survey} and by $\asymp$ in \cite{AbakumovDoubtsov15} and \cite{AbakumovDoubtsov18}. In the central statements \cite[Thm. 3.2 \& Thm. 3.3]{weightedentirecharacterization}, under the assumption that $v$ is \emph{essential} and by involving crucial results from \cite{BonetDomanskiLindstroemTaskinen98} resp. from \cite{AbakumovDoubtsov15} and \cite{AbakumovDoubtsov18}, we have shown that
\begin{equation}\label{principalequ}
v\hyperlink{ompreceq}{\preceq}w\Longleftrightarrow H^{\infty}_{v}(\CC)\subseteq H^{\infty}_{w}(\CC),
\end{equation}
with continuous inclusion and similarly if $\CC$ is replaced by $\mathbb{D}$. Note that $\Longrightarrow$ in \eqref{principalequ} is obvious, however $\Longleftarrow$ requires some additional basic assumptions on the weight $v$.

Moreover, in \cite[Sect. 3.2-3.5]{weightedentirecharacterization} we have studied the above characterization even for two types of naturally appearing \emph{weight systems:} We can either consider $t\mapsto v(ct)$, the \emph{dilatation-type,} or $t\mapsto v(t)^c$, called the \emph{exponential-type.} Here $c>0$ is a real parameter; we refer to Section \ref{weightsection} for precise definitions concerning weights and weight systems.\vspace{6pt}

The strategy has been to consider weights of the form $v=v_M: t\mapsto\exp(-\omega_M(t))$ with $M\in\RR_{>0}^{\NN}$ a given (weight) sequence such that mild standard growth assumptions are satisfied and $\omega_M$ denotes the so-called {\itshape associated weight function;} see Section \ref{assofctsect}. The advantage in this setting is that $v_M$ is based on a given sequence $M$; therefore more structure and information is involved and available. Desired properties for $v_M$ can be expressed in terms of growth and regularity assumptions, known and used in the ultradifferentiable and ultraholomorphic setting, for $M$.

Then, when an abstract weight $u$ is given, this case is reduced to the weight sequence setting by considering the so-called {\itshape associated weight sequence} $M^u$; see Section \ref{convexsection}.\vspace{6pt}

The aim of this recent work is to continue and to extend the results from \cite{weightedentirecharacterization} in several directions; however in order to conclude we follow the ideas introduced there. Roughly speaking the main goal is to replace $O$-growth relations/restrictions by their $o$-counterparts. (The spaces given by \eqref{weightedclassdef} correspond to a $O$-growth condition.) In Section \ref{equsectforo} we study the analogous characterizing results for weighted classes expressed in terms of $o$-growth restrictions; denoted by $H^{0}_v(\CC)$. Such classes have been considered also frequently in the literature; see e.g. \cite[Sect. 1.B]{BierstedtBonetTaskinen98} and \cite{BonetDomanskiLindstroemTaskinen98}. The aim is to transfer the characterizations of inclusion relations from \cite[Sect. 3]{weightedentirecharacterization} to this setting; see Theorems \ref{essentialweightcharactnew}, \ref{weightholombysequcharactsingle} and \ref{oclassesnotnew}. It turns out that the crucial characterization for (single) weight functions is precisely the same as for the weighted entire case; i.e. \eqref{principalequ} is valid for such classes as well.\vspace{6pt}

Then, in Section \ref{stronginclusionsection} we treat the situation when relation \eqref{bigOrelation} and the more general but crucial growth relations for weight systems are replaced by their corresponding $o$-growth counterparts. In this context, see the main results Theorems \ref{littleoweightsequcharact}, \ref{littleoweightsequcharactpow}, \ref{littleoweightfctcharact} and \ref{littleoweightfctcharactpow}.\vspace{6pt}

Section \ref{compdilaexposection} is dedicated to the comparison of the dilatation- and exponential-type structures. The aim is to establish a characterization when both notions coincide; see Theorems \ref{powerdilacharact} and \ref{powerdilacharactforweights}. The crucial conditions have already appeared in the weighted entire setting, see \cite[Sect. 2.8 \& 2.9]{weightedentirecharacterization} and also the literature mentioned in \cite[Sect. 2.9]{weightedentirecharacterization} concerning exponential-type systems. Moreover, these technical requirements are well-known (and standard) in the ultradifferentiable setting as well; see e.g. \cite{Komatsu73}, \cite{BraunMeiseTaylor90} and \cite{BonetMeiseMelikhov07}.\vspace{6pt}

Finally, in the Appendix \ref{erratum} we comment in detail on a gap in the proofs of the crucial results \cite[Prop. 4.7 \& 4.12]{weightedentirecharacterization} and a similar problem appearing in \cite[Prop. 4.7]{Borelmapalgebraity}. We have detected this after publishing \cite{weightedentirecharacterization} during the preparation of this article.\vspace{6pt}

\textbf{Acknowledgements.} The author thanks Prof. Jos\'{e} Bonet Solves from the Universitat Polit\`{e}cnica de Val\`{e}ncia for helpful and clarifying explanations during the preparation of this work. In particular, for forwarding the recent article \cite{Bonet2022survey}.

Moreover, the author thanks Armin Rainer from the University of Vienna for his help and joint discussions concerning Appendix \ref{erratum}.

\section{Weights, weighted spaces and growth conditions}\label{weightsection}

\subsection{Basic convention}
We use the notation $\NN=\{0,1,\dots,\}$ and $\NN_{>0}=\{1,2,\dots\}$.\vspace{6pt}

We call a weight $v$ {\itshape normalized} when $v(t)=1$ for all $t\in[0,1]$. A normalized $v$ satisfies $v(t)\le 1$ for all $t\ge 0$.

When both $v$, $w$ are (normalized) weight functions, then the product $v\cdot w$ is so, too.

For defining the class in \eqref{weightedclassdef} normalization can be assumed w.l.o.g.: Otherwise replace $v$ by a normalized weight, say $v^n$, and such that $v(t)=v^n(t)$ for all large $t>1$. This gives $H^{\infty}_v(\CC)=H^{\infty}_{v^n}(\CC)$ as l.c.v.s.

Since $w(t),v(t)\neq 0$ for all $t$, \eqref{bigOrelation} precisely means that
$$\exists\;C\ge 1\;\forall\;t\ge 0:\;\;\;w(t)\le Cv(t).$$

\subsection{Weight systems and corresponding weighted classes}\label{basicweightsection}
We revisit the notions introduced in \cite[Sect. 2.2]{weightedentirecharacterization}. Let $\underline{\mathcal{V}}=(v_n)_{n\in\NN_{>0}}$ be a non-increasing sequence of weights, i.e. $v_n\ge v_{n+1}$ for all $n$. Then we define the (LB)-space $H^{\infty}_{\underline{\mathcal{V}}}(\CC)$ of the Banach spaces $H^{\infty}_{v_n}(\CC)$, so
\begin{equation}\label{indlim}
H^{\infty}_{\underline{\mathcal{V}}}(\CC):=\varinjlim_{n\rightarrow\infty}H^{\infty}_{v_n}(\CC).
\end{equation}
Analogously, if $\overline{\mathcal{V}}=(v_n)_{n\in\NN_{>0}}$ is a non-decreasing sequence of weights, i.e. $v_n\le v_{n+1}$ for all $n$, then we define the Fr\'{e}chet-space $H^{\infty}_{\overline{\mathcal{V}}}(\CC)$ to be the projective limit of the Banach spaces $H^{\infty}_{v_n}(\CC)$, so
\begin{equation}\label{projlim}
H^{\infty}_{\overline{\mathcal{V}}}(\CC):=\varprojlim_{n\rightarrow\infty}H^{\infty}_{v_n}(\CC).
\end{equation}
Let $v$ be a weight function and $c>0$, then set
\begin{equation}\label{parameterweights}
v_c: t\mapsto v(ct),\hspace{20pt}v^c: t\mapsto v(t)^c,
\end{equation}
and write again $v$ instead of $v_1$ resp. $v^1$. $v$ is a weight function if and only if some/each weight $v_c$ resp. if and only if some/each $v^c$ is so. Note that for all $d\ge c>0$ and all large $t$ we have $v^d(t)\le v^c(t)$ since eventually $v(t)\le 1$ and we have this estimate for all $t\ge 0$ if $v$ is normalized. Moreover, $v\hyperlink{ompreceq}{\preceq}w$ if and only if $v_c\hyperlink{ompreceq}{\preceq}w_c$ if and only if $v^c\hyperlink{ompreceq}{\preceq}w^c$ for some/each $c>0$.

Consider the corresponding weight systems
\begin{equation}\label{weightsystems}
\underline{\mathcal{V}}_{\mathfrak{c}}:=(v_c)_{c\in\NN_{>0}},\hspace{20pt}\overline{\mathcal{V}}_{\mathfrak{c}}:=(v_{\frac{1}{c}})_{c\in\NN_{>0}}.
\end{equation}
Thus $\underline{\mathcal{V}}_{\mathfrak{c}}$ is a non-increasing sequence of weights, whereas $\overline{\mathcal{V}}_{\mathfrak{c}}$ is a non-decreasing sequence. Similarly, if $v\le 1$, then consider
\begin{equation}\label{weightsystems1}
\underline{\mathcal{V}}^{\mathfrak{c}}:=(v^c)_{c\in\NN_{>0}},\hspace{20pt}\overline{\mathcal{V}}^{\mathfrak{c}}:=(v^{\frac{1}{c}})_{c\in\NN_{>0}}.
\end{equation}
Since $v\le 1$ the set $\underline{\mathcal{V}}^{\mathfrak{c}}$ is a non-increasing sequence of weights, whereas $\overline{\mathcal{V}}^{\mathfrak{c}}$ is a non-decreasing sequence. The assumption $v\le 1$ holds if $v$ is normalized.

\begin{definition}
The weight systems in \eqref{weightsystems} are called to be of {\itshape dilatation-type} and the systems in \eqref{weightsystems1} of {\itshape exponential-type.}
\end{definition}

Let us write $v\hypertarget{ompreceqc}{\preceq_{\mathfrak{c}}}w$ if
\begin{equation}\label{bigOdilarelation}
\exists\;c\ge 1:\;\;\;v\hyperlink{ompreceq}{\preceq}w_c,
\end{equation}
and $v\hypertarget{simc}{\sim_{\mathfrak{c}}}w$, if
$$v\hyperlink{ompreceqc}{\preceq_{\mathfrak{c}}}w\hspace{15pt}\text{and}\hspace{15pt}w\hyperlink{ompreceqc}{\preceq_{\mathfrak{c}}}v.$$ Similarly, write $v\hypertarget{ompreceqpowc}{\preceq^{\mathfrak{c}}}w$ if
\begin{equation}\label{bigOexprelation}
\exists\;c\ge 1:\;\;\;v\hyperlink{ompreceq}{\preceq}w^c,
\end{equation}
and $v\hypertarget{simpowc}{\sim^{\mathfrak{c}}}w$, if
$$v\hyperlink{ompreceqc}{\preceq^{\mathfrak{c}}}w\hspace{15pt}\text{and}\hspace{15pt}w\hyperlink{ompreceqc}{\preceq^{\mathfrak{c}}}v.$$
Note:
\begin{itemize}
\item[$(*)$] Both $v\hyperlink{ompreceqc}{\preceq_{\mathfrak{c}}}w$ and $v\hyperlink{ompreceqpowc}{\preceq^{\mathfrak{c}}}w$ are weaker than $v\hyperlink{ompreceq}{\preceq}w$.

\item[$(*)$] $v\hyperlink{ompreceqc}{\preceq_{\mathfrak{c}}}w$ resp. $v\hyperlink{ompreceqc}{\preceq^{\mathfrak{c}}}w$ if and only if $$\exists\;c\ge 1\;\forall\;a>0:\;\;\;v_a\hyperlink{ompreceq}{\preceq}w_{ac}\;\;\;\text{resp.}\;\;\;v^a\hyperlink{ompreceq}{\preceq}w^{ac}.$$
\end{itemize}

\subsection{Weight sequences and conditions}
Let $M=(M_j)_j\in\RR_{>0}^{\NN}$. The corresponding sequence of quotients is denoted by $\mu_j:=\frac{M_j}{M_{j-1}}$, $j\ge 1$, and set $\mu_0:=1$. $M$ is called {\itshape normalized} if $1=M_0\le M_1$ and $M$ is called {\itshape log-convex} if
$$\forall\;j\in\NN_{>0}:\;M_j^2\le M_{j-1} M_{j+1},$$
equivalently if the sequence $(\mu_j)_j$ is non-decreasing. If $M_0=1$ and $M$ is log-convex, then $(M_j)^{1/j}\le\mu_j$ for all $j\in\NN_{>0}$, $j\mapsto(M_j)^{1/j}$ is non-decreasing (see e.g. \cite[Lemma 2.0.4]{diploma}) and (see e.g. \cite[Lemma 2.0.6]{diploma})
\begin{equation}\label{expsuperadd}
\forall\;j,k\in\NN:\;\;\;M_jM_k\le M_{j+k}.
\end{equation}

Now recall \cite[Def. 2.4]{weightedentirecharacterization}:

\begin{definition}\label{defweightsequ}
$M$ is called a {\itshape weight sequence} if
$$1=M_0\;\;\;\text{and}\;\;\;\lim_{j\rightarrow+\infty}(M_j)^{1/j}=+\infty.$$
\end{definition}

Moreover consider the set
$$\hypertarget{LCset}{\mathcal{LC}}:=\{M\in\RR_{>0}^{\NN}:\;M\;\text{is normalized, log-convex},\;\lim_{j\rightarrow+\infty}(M_j)^{1/j}=+\infty\},$$
i.e. $M\in\hyperlink{LCset}{\mathcal{LC}}$ if and only if $M$ is a log-convex weight sequence with $M_1\ge 1$. Note that each $M\in\hyperlink{LCset}{\mathcal{LC}}$ is non-decreasing since $\mu_j\ge 1$ for all $j$ and $j\mapsto\mu_j$ is non-decreasing.\vspace{6pt}

$M$ (with $M_0=1$) has condition {\itshape moderate growth}, denoted by \hypertarget{mg}{$(\text{mg})$}, if
$$\exists\;C\ge 1\;\forall\;j,k\in\NN:\;M_{j+k}\le C^{j+k} M_j M_k.$$
In \cite{Komatsu73} this is denoted by $(M.2)$ and also known under the name {\itshape stability under ultradifferential operators.}\vspace{6pt}

Let $M,N\in\RR_{>0}^{\NN}$, we write $M\le N$ if $M_j\le N_j$ for all $j$ and $M\hypertarget{preceq}{\preceq}N$ if $\sup_{j\in\NN_{>0}}\left(\frac{M_j}{N_j}\right)^{1/j}<+\infty$.

$M$ and $N$ are called {\itshape equivalent}, written $M\hypertarget{approx}{\approx}N$, if
$$M\hyperlink{preceq}{\preceq}N\hspace{15pt}\text{and}\hspace{15pt}N\hyperlink{preceq}{\preceq}M.$$
Property \hyperlink{mg}{$(\on{mg})$} is clearly preserved under \hyperlink{approx}{$\approx$}.

In \cite[Rem. 2.5]{weightedentirecharacterization} we have shown that for any log-convex weight sequence $M$ there exists an equivalent $N\in\hyperlink{LCset}{\mathcal{LC}}$ such that even $\mu_j=\nu_j$ for all $j$ sufficiently large.

\subsection{Associated weight function}\label{assofctsect}
Let $M\in\RR_{>0}^{\NN}$ (with $M_0=1$), then the {\itshape associated weight function} $\omega_M: \RR_{\ge 0}\rightarrow\RR\cup\{+\infty\}$ is defined by
\begin{equation*}\label{assofunc}
\omega_M(t):=\sup_{j\in\NN}\log\left(\frac{t^j}{M_j}\right)\;\;\;\text{for}\;t\in\RR_{>0},\hspace{30pt}\omega_M(0):=0.
\end{equation*}
For an abstract introduction of the associated function we refer to \cite[Chapitre I]{mandelbrojtbook}, see also \cite[Definition 3.1]{Komatsu73} and the more recent work \cite{regularnew}.

If $\liminf_{j\rightarrow+\infty}(M_j)^{1/j}>0$, then $\omega_M(t)=0$ for sufficiently small $t$, since $\log\left(\frac{t^j}{M_j}\right)<0\Leftrightarrow t<(M_j)^{1/j}$ holds for all $j\in\NN_{>0}$. (In particular, if $M_j\ge 1$ for all $j\in\NN$, then $\omega_M$ is vanishing on $[0,1]$.) Moreover, under this assumption $t\mapsto\omega_M(t)$ is a continuous non-decreasing function, which is convex in the variable $\log(t)$ and tends faster to infinity than any $\log(t^j)$, $j\ge 1$, as $t\rightarrow+\infty$. $\lim_{j\rightarrow+\infty}(M_j)^{1/j}=+\infty$ implies that $\omega_M(t)<+\infty$ for each finite $t$ which shall be considered as a basic assumption for defining $\omega_M$. In particular, this holds provided $M$ is a weight sequence.\vspace{6pt}

If $M$ is a log-convex weight sequence, then let us introduce the counting function
\begin{equation}\label{counting}
\Sigma_{M}(t):=|\{j\ge 1:\;\;\;\mu_j\le t\}|,\;\;\;t\ge 0.
\end{equation}
In this case we can compute $M$ by involving $\omega_M$ as follows, see \cite[Chapitre I, 1.4, 1.8]{mandelbrojtbook} and also \cite[Prop. 3.2]{Komatsu73}:
\begin{equation}\label{Prop32Komatsu}
M_j=\sup_{t\ge 0}\frac{t^j}{\exp(\omega_{M}(t))},\;\;\;j\in\NN.
\end{equation}
If $M$ is a weight sequence which is not log-convex, then the right-hand side of \eqref{Prop32Komatsu} yields $M^{\on{lc}}_j$ with $M^{\on{lc}}$ denoting the {\itshape log-convex minorant} of $M$ and $M^{\on{lc}}$ is a log-convex weight sequence.

Moreover, one has
\begin{equation*}\label{assovanishing}
\forall\;t\in[0,\mu_1]:\;\;\;\omega_M(t)=0,
\end{equation*}
which follows by the known integral representation formula (see \cite[1.8. III]{mandelbrojtbook} and also \cite[$(3.11)$]{Komatsu73})
\begin{equation}\label{assointrepr}
\omega_M(t)=\int_0^t\frac{\Sigma_M(u)}{u}du=\int_{\mu_1}^t\frac{\Sigma_M(u)}{u}du.
\end{equation}
In particular, if $M\in\hyperlink{LCset}{\mathcal{LC}}$, then $\omega_M$ vanishes on $[0,1]$.

Finally, for any weight sequence $M$ we have
\begin{equation}\label{lcminorantnodifference}
\omega_{M^{\on{lc}}}=\omega_{M}.
\end{equation}

\subsection{Weighted spaces of entire functions defined in terms of $\omega_M$}
For any weight sequence $M$ (in the sense of Definition \ref{defweightsequ}) and for any parameter $c>0$ we set
\begin{equation}\label{weights}
v_{M,c}(t):=\exp(-\omega_M(ct)),\;\;\;v^c_M(t):=\exp(-c\omega_M(t)),\;\;\;t\ge 0.
\end{equation}
Note that each $v_{M,c}$ and $v^c_M$ is a weight by the properties of the associated weight function. Then define
\begin{align*}
H^{\infty}_{v_{M,c}}(\CC)&:=\{f\in H(\CC): \|f\|_{v_{M,c}}:=\sup_{z\in\CC}|f(z)|v_{M,c}(|z|)<+\infty\}
\\&
=\{f\in H(\CC): \|f\|_{v_{M,c}}:=\sup_{z\in\CC}|f(z)|\exp(-\omega_M(c|z|))<+\infty\},
\end{align*}
and the class $H^{\infty}_{v^c_M}(\CC)$ is introduced accordingly. $M\in\hyperlink{LCset}{\mathcal{LC}}$ has been the standard assumption in \cite{solidassociatedweight} and in \cite[Rem. 2.6]{weightedentirecharacterization} we have seen:

\begin{itemize}
\item[$(*)$] If $M$ is a weight sequence, then in view of \eqref{lcminorantnodifference} there is no difference between the weighted classes defined in terms of $M$ or $M^{\on{lc}}$ and thus one can assume w.l.o.g. that the defining sequence is log-convex.

\item[$(*)$] Let $M$ be a log-convex weight sequence, then there exists $N\in\hyperlink{LCset}{\mathcal{LC}}$ such that as l.c.v.s.
$$\forall\;c>0:\;\;\;H^{\infty}_{v_{M,c}}(\CC)=H^{\infty}_{v_{N,c}}(\CC),\hspace{15pt}H^{\infty}_{v^c_M}(\CC)=H^{\infty}_{v^c_N}(\CC).$$
\end{itemize}

For the weight systems in \eqref{weightsystems} and \eqref{weightsystems1} we set ($M$ is fixed)
\begin{equation}\label{weightsequweightsystem}
\underline{\mathcal{M}}_{\mathfrak{c}}:=(v_{M,c})_{c\in\NN_{>0}},\hspace{10pt}\overline{\mathcal{M}}_{\mathfrak{c}}:=(v_{M,\frac{1}{c}})_{c\in\NN_{>0}},\hspace{10pt}\underline{\mathcal{M}}^{\mathfrak{c}}:=(v^c_M)_{c\in\NN_{>0}},\hspace{10pt}\overline{\mathcal{M}}^{\mathfrak{c}}:=(v^{\frac{1}{c}}_{M})_{c\in\NN_{>0}},
\end{equation}
and the (LB)-spaces from \eqref{indlim} resp. the Fr\'{e}chet-spaces from \eqref{projlim} are defined accordingly.

Now, let us recall \cite[Rem. 2.3, Lemma 2.9]{weightedentirecharacterization}:

\begin{lemma}\label{strongdilaremark}
We get the following:
\begin{itemize}
\item[$(i)$] Let $v$ be a (normalized) weight function, then
\begin{equation}\label{strongweightsrelation1}
\forall\;d>c>0:\;\;\;\lim_{t\rightarrow+\infty}\frac{v^c(t)}{v^d(t)}=\lim_{t\rightarrow+\infty}v^{c-d}(t)=+\infty.
\end{equation}

\item[$(ii)$] Let $M$ be a log-convex weight sequence. Then
\begin{equation}\label{strongweightsrelation}
\forall\;d>c>0:\;\;\;\lim_{t\rightarrow+\infty}\frac{v_{M,c}(t)}{v_{M,d}(t)}=+\infty.
\end{equation}
\end{itemize}
\end{lemma}

We close this section with the following observation:

\begin{remark}\label{BigOequivrem}
\emph{The following are equivalent by definition:}
\begin{itemize}
\item[$(*)$] \emph{$\omega_M(t)=O(\omega_N(t))$, i.e. $\omega_N\hyperlink{ompreceq}{\preceq}\omega_M$,}

\item[$(*)$] \emph{for some $c\ge 1$ we get $v_M\hyperlink{ompreceq}{\preceq}v_N^c$, i.e. $v_M\hyperlink{ompreceqpowc}{\preceq^{\mathfrak{c}}}v_N$.}
\end{itemize}
\end{remark}

\subsection{On the class of (log-)convex weights $v$ and the associated weight sequence}\label{convexsection}
We recall definitions and statements from \cite[Sect. 2.6 \& 2.7]{weightedentirecharacterization}; see also the citations there. Let $u:[0,+\infty)\rightarrow(0,+\infty)$ be a normalized (radial) weight function, then set
\begin{equation}\label{omegafromv}
\omega^u(t):=-\log(u(t)),\;\;\;t\in[0,+\infty),
\end{equation}
and we get:
\begin{itemize}
\item[$(*)$] $\omega^u:[0,+\infty)\rightarrow[0,+\infty)$ is continuous and non-decreasing,

\item[$(*)$] $\lim_{t\rightarrow+\infty}\omega^u(t)=+\infty$ and

\item[$(*)$] $\omega^u(t)=0$ for $t\in[0,1]$ (normalization).

\item[$(*)$] The fact that $u$ is rapidly decreasing is equivalent to
$$\hypertarget{om3}{(\omega_3)}:\;\;\;\log(t)=o(\omega^u(t))\hspace{15pt}t\rightarrow+\infty.$$

\item[$(*)$] We have that
\begin{equation}\label{vconvexity}
t\mapsto\omega^u(e^t)(=-\log(u(e^t)))\;\;\;\text{is convex on}\;\RR,
\end{equation}
if and only if $\varphi_{\omega^u}: t\mapsto\omega^u(e^t)$ is convex.
\end{itemize}

\hyperlink{om3}{$(\omega_3)$} is named after \cite{dissertation} and there \eqref{vconvexity} is abbreviated with $(\omega_4)$.

\begin{definition}\label{admissdef}
We call a weight $u$ {\itshape convex} if it satisfies \eqref{vconvexity}.
\end{definition}

Let $u$ be a normalized and convex weight and let $\omega^u$ be given by \eqref{omegafromv}. We introduce the sequence $M^u=(M^u_j)_{j\in\NN}$ defined by
\begin{equation}\label{vBMTweight1equ1}
M^u_j:=\sup_{t>0}\frac{t^j}{\exp(\omega^u(t))}=\sup_{t>0}t^ju(t),\;\;\;j\in\NN.
\end{equation}
Note that in order to avoid confusion instead of $v$ we denote the abstractly given weight by a different symbol than $v$, say $u$ or $w$, since we are also interested in $t\mapsto\exp(-\omega_{M^u}(t))(=v_{M^u}(t))$.\vspace{6pt}

We gather several properties:

\begin{itemize}
\item[$(i)$] $M^u\in\hyperlink{LCset}{\mathcal{LC}}$,

\item[$(ii)$] $\omega_{M^u}\hyperlink{sim}{\sim}\omega^u$, more precisely $\exists\;A\ge 1\;\forall\;t\ge 0:$
\begin{equation}\label{omegavequiv}
\frac{1}{A}(v_{M^{u}}(t))^{2}=\frac{1}{A}\exp(-2\omega_{M^u}(t))\le u(t)=\exp(-\omega^u(t))\le\exp(-\omega_{M^u}(t))=v_{M^u}(t),
\end{equation}
consequently
\begin{equation}\label{omegavequivnew}
\exists\;A\ge 1\;\forall\;c>0\;\forall\;t\ge 0:\;\;\;\frac{1}{A}v^2_{M^u,c}(t)\le u_c(t)\le v_{M^u,c}(t),\;\;\;\frac{1}{A}v^{2c}_{M^u}(t)\le u^c(t)\le v^c_{M^u}(t),
\end{equation}
which implies
\begin{equation}\label{omegavequivnewnew}
\forall\;c>0:\;\;\;H^{\infty}_{v_{M^u},c}(\CC)\subseteq H^{\infty}_{u,c}(\CC)\subseteq H^{\infty}_{v^2_{M^u},c}(\CC),\hspace{15pt}H^{\infty}_{v^c_{M^u}}(\CC)\subseteq H^{\infty}_{u^c}(\CC)\subseteq H^{\infty}_{v^{2c}_{M^u}}(\CC),
\end{equation}
with continuous inclusions. Recall that the second part in \eqref{omegavequivnew} precisely means $v_{M^u}\hyperlink{simpowc}{\sim^{\mathfrak{c}}}u$.

\item[$(iii)$] If $u\equiv v_M$ with $M\in\hyperlink{LCset}{\mathcal{LC}}$, and so $\omega^u\equiv\omega_M$, then $M^u\equiv M$.
\end{itemize}

As mentioned in \cite[Rem. 2.13]{weightedentirecharacterization} convexity for $u$ is in the whole approach only required in order to get the first estimate in \eqref{omegavequiv} resp. in \eqref{omegavequivnew}; i.e. the second inclusions in \eqref{omegavequivnewnew}.

\subsection{Technical growth conditions on weight functions}\label{modweightsection}
We revisit two (natural) but crucial growth conditions which have already appeared in the literature for both the ultradifferentiable and the weighted entire setting. For more comments in the first case we refer to \cite{BonetMeiseMelikhov07}, \cite{compositionpaper}, \cite{dissertation} and for the latter case we refer to \cite[Sect. 2.8 \& 2.9]{weightedentirecharacterization}. So let us consider

$$\hypertarget{om6}{(\omega_6)}:\;\;\; \exists\;H\ge 1\;\forall\;t\ge 0:\;\;\;2\omega(t)\le\omega(H t)+H,$$
and
$$\hypertarget{om1}{(\omega_1)}:\;\;\; \exists\;L\ge 1\;\forall\;t\ge 0:\;\;\;\omega(2t)\le L(\omega(t)+1).$$

Both conditions are named again after \cite{compositionpaper} and \cite{dissertation} and have appeared for abstractly given weight functions $\omega$ (in the sense of Braun-Meise-Taylor in \cite{BraunMeiseTaylor90}). We recall \cite[Lemmas 2.17 \& 2.20]{weightedentirecharacterization}; the proof is based on the crucial characterizations \cite[Thm. 3.1]{subaddlike} concerning \hyperlink{om1}{$(\omega_1)$} and \cite[Prop. 3.6]{Komatsu73} concerning \hyperlink{om6}{$(\omega_6)$}:

\begin{lemma}\label{techlemma}
Let $u$ be normalized and convex.
\begin{itemize}
\item[$(I)$] The following are equivalent:

\begin{itemize}
\item[$(a)$] $u$ satisfies
\begin{equation}\label{om6forv}
\exists\;H\ge 1\;\forall\;t\ge 0:\;\;\;u_H(t)=u(Ht)\le e^Hu^2(t).
\end{equation}

\item[$(b)$] The function $\omega^u$ (see \eqref{omegafromv}) satisfies \hyperlink{om6}{$(\omega_6)$}.

\item[$(c)$] The associated weight function $\omega_{M^u}$ satisfies \hyperlink{om6}{$(\omega_6)$}.

\item[$(d)$] The sequence $M^u$ satisfies \hyperlink{mg}{$(\on{mg})$}.
\end{itemize}

\item[$(II)$] The following are equivalent:

\begin{itemize}
\item[$(a)$] $u$ satisfies
\begin{equation}\label{om1forv}
\exists\;L\ge 1\;\forall\;t\ge 0:\;\;\;u^L(t)\le e^Lu_2(t).
\end{equation}

\item[$(b)$] $\omega^u$ has \hyperlink{om1}{$(\omega_1)$}.

\item[$(c)$] $\omega_{M^u}$ has \hyperlink{om1}{$(\omega_1)$}.

\item[$(d)$] $M^u$ satisfies
\begin{equation}\label{om1omegaMchar}
\exists\;L\in\NN_{>0}:\;\;\;\liminf_{j\rightarrow+\infty}\frac{(M^u_{Lj})^{1/(Lj)}}{(M^u_j)^{1/j}}>1.
\end{equation}
\end{itemize}
\end{itemize}
\end{lemma}

Of course, \eqref{om1omegaMchar} can be considered for any $M$ and we recall \cite[Def. 2.18]{weightedentirecharacterization} which is motivated by the first part of the previous result:

\begin{definition}\label{admissdef1}
A weight $v$ is called to be of \emph{moderate growth} if $v$ satisfies \eqref{om6forv}.
\end{definition}

Now $(I)$ in Lemma \ref{techlemma} can be used to transfer \eqref{strongweightsrelation} to the weight function setting.

\begin{lemma}\label{strongweightsrelationforvlem}
We have the following:
\begin{itemize}
\item[$(i)$] Let $u$ be normalized and convex. Assume that $u$ is also of moderate growth, then
\begin{equation}\label{strongweightsrelationforv}
\exists\;H\ge 1\;\forall\;c>0\;\forall\;d>Hc:\;\;\;\lim_{t\rightarrow+\infty}\frac{u_c(t)}{u_d(t)}=+\infty.
\end{equation}
Here $H$ can be chosen to be the same constant as the one appearing in \eqref{om6forv}.

\item[$(ii)$] Let $v$ and $u$ be normalized weight functions such that $v\hyperlink{ompreceqc}{\preceq_{\mathfrak{c}}}u$. If either $v$ or $u$ is in addition convex and of moderate growth, then
\begin{equation}\label{charactabstractweightcorequ}
\exists\;c\ge 1\;\exists\;H\ge 1\;\forall\;a>Hc:\;\;\;u_{a}(t)=o(v(t)),\;\;\;t\rightarrow+\infty.
\end{equation}
Here $c$ denotes the parameter from $v\hyperlink{ompreceqc}{\preceq_{\mathfrak{c}}}u$ and $H$ can be taken to be the constant appearing in \eqref{om6forv} for the particular weight.
\end{itemize}
\end{lemma}

\demo{Proof}
$(i)$ Let $d>Hc$ with $H$ the constant appearing in \eqref{om6forv}, then by combining \eqref{omegavequiv} and \eqref{om6forv} we get
$$\exists\;A\ge 1\;\exists\;H\ge 1\;\forall\;t\ge 0:\;\;\;\frac{u_c(t)}{u_d(t)}\ge\frac{1}{A}\frac{v_{{M^u},c}(t)^2}{v_{{M^u},d}(t)}\ge\frac{1}{Ae^H}\frac{v_{{M^u},Hc}(t)}{v_{{M^u},d}(t)},$$
hence \eqref{strongweightsrelation} applied to $M\equiv M^u$ yields the conclusion.\vspace{6pt}

$(ii)$ Let $c>0$ be the parameter from $v\hyperlink{ompreceqc}{\preceq_{\mathfrak{c}}}u$ and $H\ge 1$ the constant appearing in \eqref{om6forv}. If $u$ has in addition \eqref{vconvexity} and \eqref{om6forv}, then write $$\frac{u_a(t)}{v(t)}=\frac{u(at)}{u(ct)}\frac{u(ct)}{v(t)}.$$
So, by \eqref{strongweightsrelationforv} applied to $u$ we get that $\lim_{t\rightarrow+\infty}\frac{u(at)}{v(t)}=0$ for all $a>Hc$.

Similarly, if $v$ has in addition \eqref{vconvexity} and \eqref{om6forv}, then write
$$\frac{u(t)}{v(t/a)}=\frac{u(t)}{v(t/c)}\frac{v(t/c)}{v(t/a)},$$
and \eqref{strongweightsrelationforv} applied to $v$ yields $\lim_{t\rightarrow+\infty}\frac{u(at)}{v(t)}=0$ for all $a$ satisfying $c^{-1}>Ha^{-1}\Leftrightarrow a>Hc$.
\qed\enddemo

\subsection{Essential weights and optimal functions}\label{essentialsection}
In order to formulate one of the main results in the next section, i.e. Theorem \ref{essentialweightcharactnew}, and concerning crucial optimal function we have to summarize some notation. We refer to \cite[Def. 1.1, Sect. 1.B]{BierstedtBonetTaskinen98} and \cite[Sect. 3]{Bonet2022survey}; more detailed comments are also given in \cite[Sect. 3.1]{weightedentirecharacterization}.

Let a weight $v$ be given, then consider the function (radial growth condition) $w_v: \CC\rightarrow(0,+\infty)$ given by
\begin{equation}\label{corrgrowthcond}
w_v:=\frac{1}{v}
\end{equation}
and set
$$B_{w_v}(\CC):=\{f\in H(\CC): |f(z)|\le w_v(|z|),\;\forall\;z\in\CC\}.$$
The function $\widetilde{w_v}:\CC\rightarrow[0,+\infty)$ {\itshape associated with} $w_v$ is defined by
\begin{equation}\label{assoweightw}
\widetilde{w_v}(z):=\sup\{|f(z)|: f\in B_{w_v}(\CC)\},\;\;\;z\in\CC,
\end{equation}
and \cite[Observ. 1.5]{BierstedtBonetTaskinen98} gives that $\widetilde{w_v}$ is again radial. Set
\begin{equation}\label{assoweightv}
\widetilde{v}:=\frac{1}{\widetilde{w_v}}.
\end{equation}
If $v=v_M$, then we write $\widetilde{v_M}$. Similarly these notions can be introduced for arbitrary (open and connected) sets $G$. Finally, we have the following definition; see e.g. \cite[Sect. 3]{Bonet2022survey}.

\begin{definition}
A weight $v$ is called {\itshape essential} if $v\hyperlink{sim}{\sim}\widetilde{v}$.
\end{definition}

In \cite[Thm. 3.4]{weightedentirecharacterization} it has been shown that in the weighted entire case $v_M$ is essential for any log-convex weight sequence and, moreover, in this case even $v_M=\widetilde{v_M}$ is valid.

Next, we recall some optimal function belonging to weighted entire spaces; see \cite[Sect. 3.1 \& 3.2]{weightedentirecharacterization}. First, by \cite[Thm. 2]{AbakumovDoubtsov18} (see also \cite[Thm. 3]{Bonet2022survey}) we get that $u$ is essential if and only if $u$ is {\itshape approximable by the maximum of an analytic function modulus,} i.e.
\begin{equation}\label{approximable}
\exists\;f_u\in H(\CC):\;\;\;\frac{1}{u}\hyperlink{sim}{\sim} t\mapsto M(f_u,t),
\end{equation}
with
$$M(f,t):=\sup\{|f(z)|: |z|=t\};$$
see \cite[Def. 3]{AbakumovDoubtsov18} and \cite[p. 3 \& p. 5]{Bonet2022survey}. Recall that $w=\frac{1}{u}$ denotes the weight used in the notation in \cite{AbakumovDoubtsov18}. \eqref{approximable} gives that $f_u\in H^{\infty}_u(\CC)$ and more precisely this relation ensures the existence of {\itshape optimal functions} $f_u\in H^{\infty}_u(\CC)$ admitting sharp estimates from below.

On the other hand we recall \cite[Lemma 3.10, Rem. 3.12]{weightedentirecharacterization}.

\begin{lemma}\label{charholomfctlemma}
Let $M$ be a weight sequence and $c>0$. Then the function
$$\theta_{M,c}(z):=\sum_{j\ge 0}\frac{1}{2^jM_j}(cz)^j,\;\;\;z\in\CC,$$
satisfies
$$\theta_{M,c}\in H^{\infty}_{v_{M,c}}(\CC)\subseteq H^{\infty}_{\underline{\mathcal{M}}_{\mathfrak{c}}}(\CC),$$
and
\begin{equation}\label{charholomfctlemmaequ}
\forall\;c>0\;\forall\;t\ge 0:\;\;\;\exp(\omega_M(ct/2))\le|\theta_{M,c}(t)|=\theta_{M,c}(t).
\end{equation}
Moreover, for any $c\in\NN_{>0}$ let us set
\begin{equation}\label{charholomfctpow}
\theta^c_M(z):=\sum_{j\ge 0}\frac{1}{2^j(M_j)^c}z^{cj},\;\;\;z\in\CC,
\end{equation}
which satisfies $\theta^c_M\in H^{\infty}_{v^c_M}(\CC)\subseteq H^{\infty}_{\underline{\mathcal{M}}^{\mathfrak{c}}}(\CC)$ and
\begin{equation}\label{charholomfctlemmaequ1}
\forall\;c\in\NN_{>0}\;\forall\;t\ge 0:\;\;\;\exp(c\omega_M(t/2^{1/c}))\le|\theta^c_M(t)|=\theta^c_M(t).
\end{equation}
\end{lemma}

Note that $\theta_{M,c}\notin H^{\infty}_{\overline{\mathcal{M}}_{\mathfrak{c}}}(\CC)$; see \cite[Rem. 3.11]{weightedentirecharacterization}.\vspace{6pt}

Finally let $f_k(z):=z^k$, then for all $k\in\NN$ and any weight sequence $M$ we have $f_k\in H^{\infty}_{\overline{\mathcal{M}}_{\mathfrak{c}}}(\CC), H^{\infty}_{\overline{\mathcal{M}}^{\mathfrak{c}}}(\CC)$, because $|f_k(z)|v_{M,c}(|z|)=|z|^kv_M(c|z|)\rightarrow 0$ and $|f_k(z)|v^c_M(|z|)=|z|^kv^c_M(|z|)\rightarrow 0$ as $|z|\rightarrow+\infty$ for all $c>0$. For this recall that $v_M$ is rapidly decreasing. The same statement holds if $v_M$ is replaced by any abstractly given weight $v$: In fact $H_v^{\infty}(\CC)$ contains all polynomials if and only if $v$ is rapidly decreasing; see \cite[p. 3]{Bonet2022survey}.

\section{Weighted classes of entire functions with $o$-growth restriction}\label{equsectforo}
Let $v:[0,+\infty)\rightarrow(0,+\infty)$ be a given weight. We define weighted classes of entire functions with $o$-growth restriction as follows; see e.g. \cite[Sect. 1.B]{BierstedtBonetTaskinen98} and \cite{BonetDomanskiLindstroemTaskinen98}:
$$H^{0}_v(\CC):=\{f\in H(\CC): \lim_{|z|\rightarrow+\infty}|f(z)|v(|z|)=0\},$$
i.e. for all $\epsilon>0$ we can find a (radial symmetric) compact set $K\subseteq\CC$ such that $|f(z)|v(|z|)<\epsilon$ for all $z\in\CC\backslash K$. This notation has been used in \cite{Bonet2022survey} as well and similarly $H^{0}_v(\mathbb{D})$ can be introduced; see e.g. \cite{BonetDomanskiLindstroemTaskinen98}.

The spaces $H^{0}_{\underline{\mathcal{V}}_{\mathfrak{c}}}(\CC)$, $H^{0}_{\underline{\mathcal{V}}^{\mathfrak{c}}}(\CC)$, $H^{0}_{\overline{\mathcal{V}}_{\mathfrak{c}}}(\CC)$ and $H^{0}_{\overline{\mathcal{V}}^{\mathfrak{c}}}(\CC)$ are defined accordingly and similarly for the weight sequence setting. Also here normalization can be assumed w.l.o.g. for the weights and in view of \eqref{lcminorantnodifference} in the weight sequence setting log-convexity is not restricting. Concerning the growth relations we summarize:\vspace{6pt}

\begin{itemize}
	\item[$(a)$] Given two weights $u$ and $v$, then $u\hyperlink{ompreceq}{\preceq}v$ implies $H^{0}_{u}(\CC)\subseteq H^{0}_{v}(\CC)$ with continuous inclusion. Moreover, by \eqref{omegavequiv} the inclusions from \eqref{omegavequivnewnew} transfer to this setting immediately.
	
	\item[$(b)$] We also get by definition for (normalized) weight functions $u$ and $v$ the following: If $u\hyperlink{ompreceqc}{\preceq_{\mathfrak{c}}}v$, then
	$$H^{0}_{\underline{\mathcal{U}}_{\mathfrak{c}}}(\CC)\subseteq H^{0}_{\underline{\mathcal{V}}_{\mathfrak{c}}}(\CC),\hspace{20pt}H^{0}_{\overline{\mathcal{U}}_{\mathfrak{c}}}(\CC)\subseteq H^{0}_{\overline{\mathcal{V}}_{\mathfrak{c}}}(\CC),$$
	and if $u\hyperlink{ompreceqpowc}{\preceq^{\mathfrak{c}}}v$, then
	$$H^{0}_{\underline{\mathcal{U}}^{\mathfrak{c}}}(\CC)\subseteq H^{0}_{\underline{\mathcal{V}}^{\mathfrak{c}}}(\CC),\hspace{20pt}H^{0}_{\overline{\mathcal{U}}^{\mathfrak{c}}}(\CC)\subseteq H^{0}_{\overline{\mathcal{V}}^{\mathfrak{c}}}(\CC),$$
	with continuous inclusions.
	
	\item[$(c)$] Moreover, $H^{0}_{v_c}(\CC)\subseteq H^{\infty}_{v_c}(\CC)$, $H^{0}_{v^c}(\CC)\subseteq H^{\infty}_{v^c}(\CC)$ for any (normalized) weight $v$ and any $c>0$ and consequently, with continuous inclusions
	$$H^{0}_{\underline{\mathcal{V}}_{\mathfrak{c}}}(\CC)\subseteq H^{\infty}_{\underline{\mathcal{V}}_{\mathfrak{c}}}(\CC),\hspace{20pt}H^{0}_{\overline{\mathcal{V}}_{\mathfrak{c}}}(\CC)\subseteq H^{\infty}_{\overline{\mathcal{V}}_{\mathfrak{c}}}(\CC),$$
	and
	$$H^{0}_{\underline{\mathcal{V}}^{\mathfrak{c}}}(\CC)\subseteq H^{\infty}_{\underline{\mathcal{V}}^{\mathfrak{c}}}(\CC),\hspace{20pt}H^{0}_{\overline{\mathcal{V}}^{\mathfrak{c}}}(\CC)\subseteq H^{\infty}_{\overline{\mathcal{V}}^{\mathfrak{c}}}(\CC).$$
	
	\item[$(d)$] However, by the sharp growth behavior the crucial function $f_v$ from \eqref{approximable} cannot belong to the class $H^{0}_{v}(\CC)$ and so the (alternative) proof of \cite[Thm. 3.3]{weightedentirecharacterization} cannot be transferred to this setting.
	
	\item[$(e)$] On the other hand note that $f_k:=z^k\in H^{0}_v(\CC)$ for each $k\in\NN$ since $v$ is rapidly decreasing.
\end{itemize}

The crucial characterization for single weights follows, analogously as in the $H^{\infty}_v(\CC)$-case, from the more general results \cite[Prop. 2.1, Cor. 2.2]{BonetDomanskiLindstroemTaskinen98} dealing with the composition operator $C_{\varphi}(f):=f\circ\varphi$. It turns out that the characterization is provided by the same growth relation as in the $H^{\infty}_v(\CC)$-case; see \cite[Thm. 3.2]{weightedentirecharacterization}.

\begin{theorem}\label{essentialweightcharactnew}
	Let $u$ and $v$ be weights. Then the following are equivalent:
	\begin{itemize}
		\item[$(i)$] The weights are related by $\widetilde{u}\hyperlink{ompreceq}{\preceq}v$.
		
		\item[$(ii)$] The weights are related by $\widetilde{u}\hyperlink{ompreceq}{\preceq}\widetilde{v}$.
		
		\item[$(iii)$] We have $H^{0}_u(G)\subseteq H^{0}_v(G)$ (with continuous inclusion) for $G=\mathbb{D}$ and/or $G=\CC$.
	\end{itemize}
	Thus, if $u$ is essential, then the inclusion in $(iii)$ holds if and only if $u\hyperlink{ompreceq}{\preceq}v$.
\end{theorem}

When $G=\mathbb{D}$ then the relation $u\hyperlink{ompreceq}{\preceq}v$ means $v(t)\le C u(t)$ for all $0\le t<1$ and some $C\ge 1$.

\demo{Proof}
We apply \cite[Prop. 1.3. \& 2.1]{BonetDomanskiLindstroemTaskinen98} to $\varphi=\id$ and hence $C_{\id}=\id$; see also \cite[Cor. 2.2]{BonetDomanskiLindstroemTaskinen98}. Note that in \cite[Prop. 1.3 \& 2.1]{BonetDomanskiLindstroemTaskinen98} the authors deal with $G=\mathbb{D}$ but the proofs there are also valid for $G=\CC$.
\qed\enddemo

By taking into account \cite[Thm. 3.4 \& 3.5]{weightedentirecharacterization} in the weight sequence case we get the following:

\begin{theorem}\label{weightholombysequcharactsingle}
	Let $M$ and $N$ be log-convex weight sequences. Then both $v_M$ and $v_N$ are essential and the following are equivalent:
	\begin{itemize}
		\item[$(i)$] The sequences satisfy
		$$\exists\;A\ge 1\;\forall\;j\in\NN:\;\;\;N_j\le AM_j.$$
		
		\item[$(ii)$] The weights satisfy
		$$v_M\hyperlink{ompreceq}{\preceq}v_N.$$
		
		\item[$(iii)$] We have
		$$H^{0}_{v_M}(\CC)\subseteq H^{0}_{v_N}(\CC),$$
		with continuous inclusion.
	\end{itemize}
\end{theorem}

\demo{Proof}
This follows immediately by combining Theorem \ref{essentialweightcharactnew} with \cite[Thm. 3.4 \& 3.5]{weightedentirecharacterization}.
\qed\enddemo

For weight systems we find the following statement which shows that an immediate reduction to the results from \cite{weightedentirecharacterization} is possible. This is due to the fact that moving the appearing parameter allows for having ``enough space'' and so the corresponding spaces of $o$- and $O$-growth restriction coincide.

\begin{theorem}\label{oclassesnotnew}
	Let $v$ be a (normalized) weight function.
	\begin{itemize}
		\item[$(i)$] We get
		$$\forall\;d>c\ge c'>0:\;\;\;H^{\infty}_{v^{c'}}(\CC)\subseteq H^{\infty}_{v^c}(\CC)\subseteq H^{0}_{v^d}(\CC)\subseteq H^{\infty}_{v^d}(\CC),$$
		with continuous inclusions and hence
		$$H^{0}_{\underline{\mathcal{V}}^{\mathfrak{c}}}(\CC)=H^{\infty}_{\underline{\mathcal{V}}^{\mathfrak{c}}}(\CC),\hspace{15pt}H^{0}_{\overline{\mathcal{V}}^{\mathfrak{c}}}(\CC)=H^{\infty}_{\overline{\mathcal{V}}^{\mathfrak{c}}}(\CC),$$
		as locally convex vector spaces.
		
		\item[$(ii)$] If $v$ has in addition
		\begin{equation}\label{oclassesnotnewequ}
			\exists\;H\ge 1\;\forall\;c>0\;\forall\;d>Hc>0:\;\;\;\lim_{t\rightarrow+\infty}\frac{v_c(t)}{v_d(t)}=+\infty,
		\end{equation}
		then
		$$\exists\;H\ge 1\;\forall\;c\ge c'>0\;\forall\;d>Hc:\;\;\;H^{\infty}_{v_{c'}}(\CC)\subseteq H^{\infty}_{v_c}(\CC)\subseteq H^{0}_{v_d}(\CC)\subseteq H^{\infty}_{v_d}(\CC),$$
		with continuous inclusions. This implies
		$$H^{0}_{\underline{\mathcal{V}}_{\mathfrak{c}}}(\CC)=H^{\infty}_{\underline{\mathcal{V}}_{\mathfrak{c}}}(\CC),\hspace{15pt}H^{0}_{\overline{\mathcal{V}}_{\mathfrak{c}}}(\CC)=H^{\infty}_{\overline{\mathcal{V}}_{\mathfrak{c}}}(\CC),$$
		as locally convex vector spaces.
	\end{itemize}
\end{theorem}

\demo{Proof}
$(i)$ The first and the third inclusion is trivial, the second one follows from \eqref{strongweightsrelation1}.

$(ii)$  The first and the third inclusion is trivial, the second one follows from \eqref{oclassesnotnewequ} with the same $H$.
\qed\enddemo

We comment in detail on the consequences of this result:

\begin{itemize}
	\item[$(I)$] By Lemma \ref{strongdilaremark}, see \eqref{strongweightsrelation}, requirement \eqref{oclassesnotnewequ} is satisfied with $H=1$ in the weight sequence setting when $M$ is a given log-convex weight sequence and so as l.c.v.s. $$H^{0}_{\underline{\mathcal{M}}_{\mathfrak{c}}}(\CC)=H^{\infty}_{\underline{\mathcal{M}}_{\mathfrak{c}}}(\CC),\hspace{15pt}H^{0}_{\overline{\mathcal{M}}_{\mathfrak{c}}}(\CC)=H^{\infty}_{\overline{\mathcal{M}}_{\mathfrak{c}}}(\CC).$$
	Summarizing, in the weight sequence setting the dilatation- and exponential-type spaces with $o$- and $O$-growth conditions coincide.
	
	\item[$(II)$] Recall that \eqref{oclassesnotnewequ} is not clear for abstractly given weights $v$. However, by $(i)$ in Lemma \ref{strongweightsrelationforvlem} this property holds provided $v$ is a normalized and convex weight function of \emph{moderate growth.} Here we can choose $H$ to be the constant appearing in \eqref{om6forv} (see \eqref{strongweightsrelationforv}) and note that in this case $H>1$; see also \cite[Sect. 2.8]{weightedentirecharacterization}.

\item[$(III)$] Let $u$ be a normalized (and convex) weight. Then, by taking into account the inclusions from Theorem \ref{oclassesnotnew}, comment $(I)$ before, the facts that $v_{M^u}$ is essential, see \cite[Thm. 3.4 $(i)$]{weightedentirecharacterization}, and that \eqref{omegavequivnewnew} transfers to the $o$-setting, we get for the optimal functions from \eqref{approximable}
    $$\forall\;c>1:\;\;\;f_{v_{M^u}}\in H^{\infty}_{v_{M^u}}(\CC)\subseteq H^{0}_{v_{M^u}^c}(\CC)\subseteq H^{0}_{u^c}(\CC),$$
    and
    $$\forall\;c>1:\;\;\;f_{v_{M^u}}\in H^{\infty}_{v_{M^u}}(\CC)\subseteq H^{0}_{v_{M^u},c}(\CC)\subseteq H^{0}_{u,c}(\CC).$$
    This should be compared with comment $(d)$ above in this section.

    Moreover, if $M$ is a weight sequence, then for the functions from Lemma \ref{charholomfctlemma}, by the inclusions stated in Theorem \ref{oclassesnotnew} and comment $(I)$, we have that
\begin{equation}\label{strongthetaequ}
\forall\;c'>c>0:\;\;\;\theta_{M,c}\in H^{\infty}_{v_{M,c}}(\CC)\subseteq H^{0}_{v_{M,c'}}(\CC)\subseteq H^{0}_{\underline{\mathcal{M}}_{\mathfrak{c}}}(\CC),
\end{equation}
and
\begin{equation}\label{strongthetaequ1}
\forall\;c\in\NN_{>0}\;\forall\;c'>c:\;\;\;\theta^c_M\in H^{\infty}_{v^c_M}(\CC)\subseteq H^{0}_{v^{c'}_M}(\CC)\subseteq H^{\infty}_{\underline{\mathcal{M}}^{\mathfrak{c}}}(\CC).
\end{equation}
\end{itemize}

Now, by using Theorem \ref{oclassesnotnew} and these additional observations we see that the characterizing results concerning the inclusion relations for $O$-type-classes from \cite[Sect. 3.2-3.5]{weightedentirecharacterization}, i.e. \cite[Thm. 3.15, 3.20, 3.25 \& 3.28]{weightedentirecharacterization}, transfer to the corresponding analogues for classes with $o$-growth restriction.\vspace{6pt}

More precisely, for the exponential-type systems this is immediate by $(i)$ in Theorem \ref{oclassesnotnew} and by $(ii)$ in this result we have this also for the dilatation-type weight sequence setting. Concerning the abstract dilatation-type weight function setting we point out that, when assuming that \emph{both} weights are convex and of moderate growth, then $(ii)$ in Theorem \ref{oclassesnotnew} directly yields the conclusion. However, we show now that the weaker assumptions from \cite[Thm. 3.25]{weightedentirecharacterization}, i.e. requiring moderate growth for one of the appearing weights, are sufficient to conclude.\vspace{6pt}

Let $u$ and $w$ be normalized weights such that $u$ is convex and that either $u$ or $w$ is of moderate growth. We distinguish; similar techniques are also required when transferring the proofs from the $O$- to the $o$-setting:
		
\begin{itemize}
	\item[$(\ast)$] \emph{The inductive case:} Assume now that $H^{0}_{\underline{\mathcal{U}}_{\mathfrak{c}}}(\CC)\subseteq H^{0}_{\underline{\mathcal{W}}_{\mathfrak{c}}}(\CC)$ (as sets). Since the inclusions in \eqref{omegavequivnewnew} transfer to the $o$-setting we have
	$$H^{0}_{v_{M^u}}(\CC)\subseteq H^{0}_{u}(\CC)\subseteq H^{0}_{\underline{\mathcal{U}}_{\mathfrak{c}}}(\CC)\subseteq H^{0}_{\underline{\mathcal{W}}_{\mathfrak{c}}}(\CC).$$
	Then consider $\theta_{M,c}$ for some $0<c<1$ arbitrary but from now on fixed. In view of \eqref{strongthetaequ} the previous inclusions yield the estimate
	$$\exists\;d>0\;\exists\;D\ge 1\;\forall\;t\ge 0:\;\;\;\exp(\omega_{M^u}(ct/2))\le|\theta_{M^u,c}(t)|=\theta_{M^u,c}(t)\le D(w(dt))^{-1}.$$
	Using this and involving \eqref{Prop32Komatsu} for $M^u$ and \eqref{vBMTweight1equ1} for $M^w$ we obtain for all $j\in\NN$:
	$$M^u_j=\sup_{t\ge 0}\frac{(ct/2)^j}{\exp(\omega_{M^u}(ct/2))}\ge (c/2)^jD^{-1}\sup_{t\ge 0}t^jw(dt)=\left(\frac{c}{2d}\right)^jD^{-1}M^w_j;$$
	i.e. $M^w\le D(2d/c)^jM^u_j$ for all $j\in\NN$ is verified. By definition of associated weight functions this yields $\omega_{M^u}(t)\le\omega_{M^w}(2dt/c)+\log(D)$ and so $v_{M^w}(2dt/c)\le Dv_{M^u}(t)$ for all $t\ge 0$. Finally, we apply the first part of \eqref{omegavequiv} to $u$, and for which \emph{convexity} is required, and the second part there to $w$ in order to get
	$$\exists\;d>0\;\exists\;D,B\ge 1\;\forall\;t\ge 0:\;\;\;w(2dt/c)\le v_{M^w}(2dt/c)\le Dv_{M^u}(t)\le DB\sqrt{u(t)}.$$
	Thus a variant of \cite[$(3.36),(3.39)$]{weightedentirecharacterization} is shown and by applying moderate growth to one of the weights (the proof of) \cite[Prop. 3.22 $(ii)$]{weightedentirecharacterization} yields $u\hyperlink{ompreceqc}{\preceq_{\mathfrak{c}}}w$.

\item[$(\ast)$] {\itshape The projective case:} Note that all monomials $f_k(z):=z^k$ are also contained in $H^{0}_{\overline{\mathcal{V}}_{\mathfrak{c}}}(\CC)$, so the proof of \cite[Prop. 3.24]{weightedentirecharacterization} can be repeated and hence yields $u\hyperlink{ompreceqc}{\preceq_{\mathfrak{c}}}w$, too.
\end{itemize}

We close this section with the following observation:

\begin{remark}
\emph{The main results from \cite[Sect. 4]{weightedentirecharacterization} concerning stability under the action of the bilinear point-wise multiplication operator $\mathfrak{m}: (f,g)\mapsto f\cdot g$ transfer to the $o$-setting.}

\emph{More precisely, by $(i)$ in Theorem \ref{oclassesnotnew} the exponential-type result \cite[Prop. 4.1]{weightedentirecharacterization} is immediate. Concerning the dilatation-type results \cite[Thm. 4.8 \& 4.13]{weightedentirecharacterization} for the weight sequence setting also nothing has to be shown since both notions coincide as seen before. Finally, concerning the dilatation-type weight function setting, by taking into account \eqref{strongthetaequ} and following the ideas in the inductive case before, in the proof of \cite[Prop. 4.11]{weightedentirecharacterization} we have to consider the functions $\theta_{M^r,c'}$, $\theta_{M^u,c'}$ with $0<c'<\min\{c_1,c_2\}$ and $c_1$, $c_2$ being the parameters appearing in the $o$-growth analogue assumption in \cite[Prop. 4.11 $(i)$]{weightedentirecharacterization}. Since the monomials are contained in the spaces, the proof of \cite[Prop. 4.12]{weightedentirecharacterization} directly transfers to the $o$-setting and so also here the $o$-analogue of \cite[Thm. 4.13]{weightedentirecharacterization} is valid.}

\emph{However, if the characterization \cite[Thm. 4.13]{weightedentirecharacterization} established in the dilatation-type weight function setting for the weight $u$ holds, i.e. $u$ satisfies \eqref{om6forv}, then by Lemma \ref{strongweightsrelationforvlem} we see that in this case $H^{0}_{\underline{\mathcal{U}}_{\mathfrak{c}}}(\CC)$ and $H^{0}_{\overline{\mathcal{U}}_{\mathfrak{c}}}(\CC)$ are coinciding with their $O$-counterparts as l.c.v.s., too.}

\emph{Finally, note that also the proof of \cite[Cor. 4.14]{weightedentirecharacterization} applies to the $o$-growth setting. According to \cite[Prop. 4.12, Cor. 4.14]{weightedentirecharacterization}, please see also the comments and corrections given in Appendix \ref{erratum}.}
\end{remark}

\section{Strong inclusion relations for weighted spaces of entire functions}\label{stronginclusionsection}
The aim of this section is to replace the growth relations from Section \ref{basicweightsection} by their stronger $o$-growth counterparts and prove analogous versions of the main characterizing results from \cite[Sect. 3.2-3.5]{weightedentirecharacterization}.

\subsection{Strong growth relations for weight functions and weight sequences}
Let $v$ and $w$ be (normalized) weight functions. We write $v\hypertarget{omtriangle}{\vartriangleleft}w$ if
$$w(t)=o(v(t))\;\;\;\text{as}\;\;\;t\rightarrow+\infty.$$
Obviously $v\hyperlink{omtriangle}{\vartriangleleft}w$ implies $v\hyperlink{ompreceq}{\preceq}w$ and \hyperlink{omtriangle}{$\vartriangleleft$} is neither reflexive nor symmetric. Next we introduce
$$v\hypertarget{trianglec}{\vartriangleleft_{\mathfrak{c}}}w\;\;\;\text{resp.}\;\;\;v\hypertarget{trianglecpow}{\vartriangleleft^{\mathfrak{c}}}w,$$ defined by
$$\forall\;c>0:\;\;\;v\hyperlink{omtriangle}{\vartriangleleft}w_c\;\;\;\text{resp.}\;\;\;v\hyperlink{omtriangle}{\vartriangleleft}w^c.$$

\begin{remark}\label{strongconfusionrem}
	\emph{We gather some facts:}
	\begin{itemize}
		\item[$(i)$] \emph{In both \hyperlink{trianglec}{$\vartriangleleft_{\mathfrak{c}}$} and \hyperlink{trianglecpow}{$\vartriangleleft^{\mathfrak{c}}$} we are interested in small $c(\le 1)$.}
		
		\item[$(ii)$] \emph{Clearly, $v\hyperlink{trianglec}{\vartriangleleft_{\mathfrak{c}}}w$ resp. $v\hyperlink{trianglecpow}{\vartriangleleft^{\mathfrak{c}}}w$ is equivalent to}
		\begin{equation}\label{strongweightrelation}
			\forall\;c,d>0:\;\;\;v_d\hyperlink{mtriangle}{\vartriangleleft}w_c\;\;\;\text{resp.}\;\;\;v^d\hyperlink{mtriangle}{\vartriangleleft}w^c.
		\end{equation}
		
		\item[$(iii)$] \emph{Since by \eqref{strongweightsrelation1} we have $w^c\hyperlink{omtriangle}{\vartriangleleft}w^d$ for all $0<c<d$, condition $v\hyperlink{trianglecpow}{\vartriangleleft^{\mathfrak{c}}}w$ is equivalent to}
		\begin{equation}\label{strongweightrelationpow}
			\forall\;c>0:\;\;\;v\hyperlink{ompreceq}{\preceq}w^c.
		\end{equation}
	\end{itemize}
\end{remark}

On the other hand, let us write $M\hypertarget{mtriangle}{\vartriangleleft}N$ if
$$\forall\;h>0\;\exists\;C_h\ge 1\;\forall\;j\in\NN:\; M_j\le C_hh^jN_j\Longleftrightarrow\lim_{j\rightarrow+\infty}\left(\frac{M_j}{N_j}\right)^{1/j}=0.$$
Obviously $M\hyperlink{mtriangle}{\vartriangleleft}N$ implies $M\hyperlink{preceq}{\preceq}N$ and \hyperlink{mtriangle}{$\vartriangleleft$} is neither reflexive nor symmetric.

We gather consequences of these relations; compare this also with Lemma \ref{strongweightsrelationforvlem}.

\begin{lemma}\label{charactabstractweightcor}
	Let $M,N$ be weight sequences and assume that $M$ is log-convex. Then the following are equivalent:
	\begin{itemize}
		\item[$(a)$] $M\hyperlink{mtriangle}{\vartriangleleft}N$ is valid,
		
		\item[$(b)$] $v_N\hyperlink{trianglec}{\vartriangleleft_{\mathfrak{c}}}v_M$ is valid,
		
		\item[$(c)$] $v_N\hyperlink{ompreceq}{\preceq}v_{M,c}$ is valid for all $c>0$.
	\end{itemize}
\end{lemma}

\demo{Proof}
$(a)\Rightarrow(b)$ $M\hyperlink{mtriangle}{\vartriangleleft}N$ implies by definition
\begin{equation}\label{charactabstractweightcorequ}
	\forall\;h>0\;\exists\;C_h\ge 1\;\forall\;t\ge 0:\;\;\;\omega_N(t)\le\omega_M(ht)+\log(C_h).
\end{equation}
Let now $c>c'>0$ be given and fixed; the case $c'\ge c$ follows then. Next we choose $d>0$ and $h>0$ such that $c'>d$ and $hc\le d$. So $\omega_N(ct)\le\omega_M(dt)+\log(C_h)$ for all $t\ge 0$ follows and \eqref{assointrepr} yields
$$\omega_M(c't)-\omega_M(dt)=\int_{dt}^{c't}\frac{\Sigma_M(u)}{u}du\ge\Sigma_M(dt)\int_{dt}^{c't}\frac{1}{u}du=\Sigma_M(dt)\log(c'/d)\rightarrow+\infty\;\;\;\text{as}\;t\rightarrow+\infty,$$
by the definition of $\Sigma_M$ in \eqref{counting}. Thus
$$\forall\;\epsilon>0\;\exists\;t_{\epsilon}>0\;\forall\;t\ge t_{\epsilon}:\;\;\;\omega_M(dt)+\frac{1}{\epsilon}\le\omega_M(c't),$$
and gathering everything we get
$$\forall\;\epsilon>0\;\exists\;t_{\epsilon}>0\;\forall\;t\ge t_{\epsilon}:\;\;\;\omega_N(ct)\le\omega_M(dt)+\log(C_h)\le\log(C_h)-\frac{1}{\epsilon}+\omega_M(c't).$$
This verifies the first part of \eqref{strongweightrelation} since $C_h$ is only depending on given (fixed) choices $c,c'$ but not on $\epsilon\rightarrow 0$.\vspace{6pt}

$(b)\Rightarrow(c)$ This is trivial.\vspace{6pt}

$(c)\Rightarrow(a)$ If $v_N\hyperlink{ompreceq}{\preceq}v_{M,c}$ is valid for all $c>0$, then
$$\forall\;c>0\;\exists\;C\ge 1\;\forall\;t\ge 0:\;\;\; v_{M,c}(t)\le Cv_N(t)\Leftrightarrow\frac{1}{\exp(\omega_M(ct))}\le C\frac{1}{\exp(\omega_N(t))}.$$
Thus \eqref{Prop32Komatsu} yields for all $j\in\NN$
\begin{align*}
	M_j&=\sup_{t\ge 0}\frac{t^j}{\exp(\omega_M(t))}=c^j\sup_{t\ge 0}\frac{t^j}{\exp(\omega_M(ct))}\le Cc^j\sup_{t\ge 0}\frac{t^j}{\exp(\omega_N(t))}=Cc^jN^{\on{lc}}_j\le Cc^jN_j,
\end{align*}
which gives $M\hyperlink{mtriangle}{\vartriangleleft}N$.
\qed\enddemo

\begin{lemma}\label{charactabstractweightcor1}
	Let $M,N\in\hyperlink{LCset}{\mathcal{LC}}$.
	\begin{itemize}
		\item[$(a)$] If either $M$ or $N$ has in addition \hyperlink{mg}{$(\on{mg})$}, then $v_N\hyperlink{trianglec}{\vartriangleleft_{\mathfrak{c}}}v_M$ implies $v_N\hyperlink{trianglecpow}{\vartriangleleft^{\mathfrak{c}}}v_M$.
		
		\item[$(b)$] If either $M$ or $N$ has in addition \eqref{om1omegaMchar}, then $v_N\hyperlink{trianglecpow}{\vartriangleleft^{\mathfrak{c}}}v_M$ implies $v_N\hyperlink{trianglec}{\vartriangleleft_{\mathfrak{c}}}v_M$.
		
		\item[$(c)$] The following are equivalent:
		\begin{itemize}
			\item[$(i)$] $v_N\hyperlink{trianglecpow}{\vartriangleleft^{\mathfrak{c}}}v_M$,
			
			\item[$(ii)$] $\omega_N(t)=o(\omega_M(t))$ as $t\rightarrow+\infty$,
			
			\item[$(iii)$] we have that
			\begin{equation}\label{tildestrong}
				\forall\;c\in\NN_{>0}\;\exists\;A\ge 1\;\forall\;j\in\NN:\;\;\;\widetilde{M}^c_j:=(M_{cj})^{1/c}\le AN_j.
			\end{equation}
		\end{itemize}
	\end{itemize}
\end{lemma}

\demo{Proof}
The proof of $(a)$ and $(b)$ should be compared with the techniques in Proposition \ref{secondweighfctprop}.

$(a)$ If $M$ has \hyperlink{mg}{$(\on{mg})$} then we get \hyperlink{om6}{$(\omega_6)$} for $\omega_M$; recall \cite[Prop. 3.6]{Komatsu73}. Let $1>c>0$ be given, then iterating \hyperlink{om6}{$(\omega_6)$} yields
$$\exists\;H\ge 1\;\forall\;t\ge 0:\;\;\;\omega_M(t)\le c\omega_M(Ht)+cH.$$
By Lemma \ref{charactabstractweightcor} relation $v_N\hyperlink{trianglec}{\vartriangleleft_{\mathfrak{c}}}v_M$ implies $M\hyperlink{mtriangle}{\vartriangleleft}N$ and so \eqref{charactabstractweightcorequ}. We apply this to $h:=H^{-1}$ and get
$$\exists\;C_{H^{-1}}\ge 1\;\forall\;t\ge 0:\;\;\;\omega_N(t)\le\omega_M(H^{-1}t)+\log(C_{H^{-1}})\le c\omega_M(t)+cH+\log(C_{H^{-1}});$$
i.e. $v^c_M(t)\le e^{cH}C_{H^{-1}}v_N(t)$. Thus \eqref{strongweightrelationpow} is shown and in view of $(iii)$ in Remark \ref{strongconfusionrem} we are done.

On the other hand, if $N$ has \hyperlink{mg}{$(\on{mg})$}, then similarly by iterating \hyperlink{om6}{$(\omega_6)$} for $\omega_N$ and \eqref{charactabstractweightcorequ} we have
$$\omega_N(t)\le c\omega_N(Ht)+cH\le c\omega_M(Hht)+c\log(C_h)+cH.$$
\eqref{charactabstractweightcorequ} again applied to $h:=H^{-1}$ yields the conclusion.\vspace{6pt}

$(b)$ If $M$ satisfies \eqref{om1omegaMchar}, then $\omega_M$ has \hyperlink{om1}{$(\omega_1)$}; recall \cite[Thm. 3.1]{subaddlike}. Let $1>h>0$ be given, then iterating \hyperlink{om1}{$(\omega_1)$} gives
$$\exists\;L\ge 1\;\forall\;t\ge 0:\;\;\;\omega_M(t)\le L\omega_M(th)+L.$$
$v_N\hyperlink{trianglecpow}{\vartriangleleft^{\mathfrak{c}}}v_M$ is equivalent to \eqref{strongweightrelationpow}; see $(iii)$ in Remark \ref{strongconfusionrem}. We apply this relation to $c:=L^{-1}$ and get
$$\exists\;C\ge 1\;\forall\;t\ge 0:\;\;\;\omega_N(t)\le c\omega_M(t)+C\le\omega_M(th)+1+C\Longrightarrow v_{M,h}(t)\le e^{1+C}v_N(t).$$
Since $h>0$ was arbitrary we have verified $v_N\hyperlink{ompreceq}{\preceq}v_{M,h}$ for all $h>0$ (small) and thus by Lemma \ref{charactabstractweightcor} we are done.

On the other hand, if $N$ satisfies \eqref{om1omegaMchar}, then similarly $\omega_N(t)\le L\omega_N(th)+L$ and so again \eqref{strongweightrelationpow} applied to $c:=L^{-1}$ gives
$$\omega_N(t)\le L\omega_N(th)+L\le Lc\omega_M(th)+LC+L\Longrightarrow v_{M,h}(t)\le e^{L(C+1)}v_N(t).$$

$(c)(i)\Leftrightarrow(ii)$ First, analogously to
Remark \ref{BigOequivrem} relation $\omega_N(t)=o(\omega_M(t))$ as $t\rightarrow+\infty$, i.e. $\omega_M\hyperlink{omtriangle}{\vartriangleleft}\omega_N$, holds if and only if for all $0<c\le 1$ we get $v_N\hyperlink{ompreceq}{\preceq}v_M^c$; i.e. \eqref{strongweightrelationpow}. By $(ii),(iii)$ in Remark \ref{strongconfusionrem} this is equivalent to $v_N\hyperlink{omtriangle}{\vartriangleleft}v_M^c$ for all $0<c\le 1$ and so precisely $v_N\hyperlink{trianglecpow}{\vartriangleleft^{\mathfrak{c}}}v_M$.\vspace{6pt}

$(c)(ii)\Leftrightarrow(iii)$ This is the second part of \cite[Lemma 6.5]{PTTvsmatrix}.
\qed\enddemo

\subsection{The dilatation-type weight sequence case}\label{strongweightsequsect}

We start with the following result.

\begin{proposition}\label{charactprop2foro}
	Let $M,N$ be weight sequences and assume that $M$ is log-convex. Consider the following assertions:
	\begin{itemize}
		\item[$(i)$] The weights satisfy $M\hyperlink{mtriangle}{\vartriangleleft}N$.
		
		\item[$(ii)$] The weights satisfy $v_N\hyperlink{trianglec}{\vartriangleleft_{\mathfrak{c}}}v_M$.
		
		\item[$(iii)$] We have that
		\begin{equation}\label{charactprop2equforo}
			\forall\;c,c'>0:\;\;\;H^{\infty}_{v_{N,c}}(\CC)\subseteq H^{\infty}_{v_{M,c'}}(\CC),
		\end{equation}
		with continuous inclusion.
		\item[$(iv)$] We have that
		\begin{equation}\label{charactprop2equ1foro}
			H^{\infty}_{\underline{\mathcal{N}}_{\mathfrak{c}}}(\CC)\subseteq H^{\infty}_{\overline{\mathcal{M}}_{\mathfrak{c}}}(\CC),
		\end{equation}
		with continuous inclusion.
	\end{itemize}
	Then $(i)\Leftrightarrow(ii)\Rightarrow(iii)\Rightarrow(iv)$ holds.
\end{proposition}

\demo{Proof}
$(i)\Leftrightarrow(ii)$ is $(a)\Leftrightarrow(b)$ in Lemma \ref{charactabstractweightcor} and $(ii)\Rightarrow(iii)\Rightarrow(iv)$ is clear.
\qed\enddemo

Now we treat the converse statement.

\begin{proposition}\label{charactprop3foro}
	Let $M,N$ be weight sequences and assume that $M$ is log-convex. If we have
	\begin{equation}\label{charactprop3equforo}
		\exists\;c>0\;\forall\;c'>0:\;\;\;H^{\infty}_{v_{N,c}}(\CC)\subseteq H^{\infty}_{v_{M,c'}}(\CC),
	\end{equation}
	and if we assume that the inclusion holds as sets, then $M\hyperlink{mtriangle}{\vartriangleleft}N$ is valid.
\end{proposition}
Note that \eqref{charactprop2equforo} is stronger than \eqref{charactprop3equforo}.

\demo{Proof}
Let $c>0$ be the parameter and by assumption $\theta_{N,c}\in H^{\infty}_{v_{M,c'}}(\CC)$ for any $c'>0$. Hence
$$\forall\;c'>0\;\exists\;C\ge 1\;\forall\;t\ge 0:\;\;\;\exp(\omega_N(ct/2))\le|\theta_{N,c}(t)|\le C\exp(\omega_M(c't)).$$
By using \eqref{Prop32Komatsu} we obtain for all $j\in\NN$ that
\begin{align*}
	M_j&=\sup_{t\ge 0}\frac{t^j}{\exp(\omega_{M}(t))}=\sup_{s\ge 0}\frac{(c's)^j}{\exp(\omega_{M}(c's))}\le C(c')^j\sup_{s\ge 0}\frac{s^j}{\exp(\omega_{N}(cs/2))}
	\\&
	=C(2c'/c)^j\sup_{u\ge 0}\frac{u^j}{\exp(\omega_{N}(u))}=C(2c'/c)^jN^{\on{lc}}_j\le C(2c'/c)^jN_j.
\end{align*}
Then let $c'\rightarrow 0$ and so $M\hyperlink{mtriangle}{\vartriangleleft}N$ is verified; for this note that $c$ is given and fixed.\vspace{6pt}

A second and independent argument, see also \cite[Prop. 3.13]{weightedentirecharacterization}, is to apply \cite[Thm. 3.5]{weightedentirecharacterization} to the inclusion \eqref{charactprop3equforo} and to the weights $v_{N^{\frac{1}{c}}}$, $v_{M^{\frac{1}{c'}}}$ with $N^{\frac{1}{c}}_j:=c^{-j}N_j$, see \cite[Rem. 2.10 $(a)$]{weightedentirecharacterization}. This yields
$$\exists\;c>0\;\forall\;c'>0\;\exists\;A\ge 1\;\forall\;j\in\NN:\;\;\;\frac{1}{(c')^j}M_j=M^{\frac{1}{c'}}_j\le AN^{\frac{1}{c}}_j=A\frac{1}{c^j}N_j.$$
As $c'\rightarrow 0$ this verifies again $M\hyperlink{mtriangle}{\vartriangleleft}N$. On the one hand, here the technical constant $2$ from \eqref{charholomfctlemmaequ} disappears but which does not effect the conclusion. On the other hand, in order to apply \cite[Thm. 3.5]{weightedentirecharacterization} we also require log-convexity for $N$ but recall that in view of \eqref{lcminorantnodifference} the failure of log-convexity is negligible.
\qed\enddemo

Gathering Lemma \ref{charactabstractweightcor} and Propositions \ref{charactprop2foro} and \ref{charactprop3foro} we arrive at the following characterization:

\begin{theorem}\label{littleoweightsequcharact}
	Let $M,N$ be weight sequences and assume that $M$ is log-convex. Then the following are equivalent:
	\begin{itemize}
		\item[$(i)$] We have $M\hyperlink{mtriangle}{\vartriangleleft}N$.
		
		\item[$(ii)$] We have $v_N\hyperlink{trianglec}{\vartriangleleft_{\mathfrak{c}}}v_M$.
		
		\item[$(iii)$] We have $v_N\hyperlink{ompreceq}{\preceq}v_{M,c}$ for all $c>0$.
		
		\item[$(iv)$] \eqref{charactprop2equforo} holds; i.e.
		$$\forall\;c,c'>0:\;\;\;H^{\infty}_{v_{N,c}}(\CC)\subseteq H^{\infty}_{v_{M,c'}}(\CC),$$
		with continuous inclusion.
		
		\item[$(v)$] \eqref{charactprop3equforo} holds; i.e.
		$$\exists\;c>0\;\forall\;c'>0:\;\;\;H^{\infty}_{v_{N,c}}(\CC)\subseteq H^{\infty}_{v_{M,c'}}(\CC),$$
		with continuous inclusion.
		
		\item[$(vi)$] \eqref{charactprop2equ1foro} holds; i.e.
		$$H^{\infty}_{\underline{\mathcal{N}}_{\mathfrak{c}}}(\CC)\subseteq H^{\infty}_{\overline{\mathcal{M}}_{\mathfrak{c}}}(\CC),$$
		with continuous inclusion.
	\end{itemize}
\end{theorem}

\subsection{The exponential-type weight sequence case}
We start with the following result.

\begin{proposition}\label{charactprop2foropow}
	Let $M,N\in\hyperlink{LCset}{\mathcal{LC}}$ be given and consider the following assertions:
	\begin{itemize}
		\item[$(i)$] $M\hyperlink{mtriangle}{\vartriangleleft}N$ is valid and assume that either $M$ or $N$ has in addition \hyperlink{mg}{$(\on{mg})$}.

\item[$(ii)$] $M$ and $N$ are related by \eqref{tildestrong}.
		
		\item[$(iii)$] We have $v_N\hyperlink{trianglecpow}{\vartriangleleft^{\mathfrak{c}}}v_M$.
		
		\item[$(iv)$] We have
		\begin{equation}\label{charactprop2equforopow}
			\forall\;c,c'>0:\;\;\;H^{\infty}_{v^c_{N}}(\CC)\subseteq H^{\infty}_{v^{c'}_{M}}(\CC),
		\end{equation}
		with continuous inclusion.
		
		\item[$(v)$] We have
		\begin{equation}\label{charactprop2equ1foropow}
			H^{\infty}_{\underline{\mathcal{N}}^{\mathfrak{c}}}(\CC)\subseteq H^{\infty}_{\overline{\mathcal{M}}^{\mathfrak{c}}}(\CC),
		\end{equation}
		with continuous inclusion.
	\end{itemize}
	Then $(i)\Rightarrow(ii)\Leftrightarrow(iii)\Rightarrow(iv)\Rightarrow(v)$ holds.
\end{proposition}

\demo{Proof}
$(i)\Rightarrow(iii)$ follows by combining Lemma \ref{charactabstractweightcor} and $(a)$ in Lemma \ref{charactabstractweightcor1}, $(ii)\Leftrightarrow(iii)$ holds by $(c)$ in Lemma \ref{charactabstractweightcor1}. $(iii)\Rightarrow(iv)$ follows immediately by definition of relation $v_N\hyperlink{trianglecpow}{\vartriangleleft^{\mathfrak{c}}}v_M$ and $(iv)\Rightarrow(v)$ is clear.
\qed\enddemo

\begin{proposition}\label{charactprop3foropow}
	Let $M,N\in\hyperlink{LCset}{\mathcal{LC}}$ be given. Then the following are equivalent:
	\begin{itemize}
		\item[$(i)$]  We have (as sets)
		\begin{equation*}\label{charactprop3equ0foropow}
			\forall\;c\in\NN_{>0}\;\forall\;c'>0:\;\;\;H^{\infty}_{v^c_N}(\CC)\subseteq H^{\infty}_{v^{c'}_M}(\CC).
		\end{equation*}
		
		\item[$(ii)$] We have (as sets)
		\begin{equation}\label{charactprop3equforopow}
			\exists\;c\in\NN_{>0}\;\forall\;c'>0:\;\;\;H^{\infty}_{v^c_N}(\CC)\subseteq H^{\infty}_{v^{c'}_M}(\CC).
		\end{equation}
		
		\item[$(iii)$] The weight sequences satisfy
		\begin{equation}\label{charactprop3equ1foropow}
			\exists\;c\in\NN_{>0}\;\forall\;d\in\NN_{>0}\;\exists\;C\ge 1\;\forall\;j\in\NN:\;\;\;\widetilde{M}^{cd}_j=(M_{cdj})^{1/(cd)}\le CN_j.
		\end{equation}
		\item[$(iv)$] The weight sequences satisfy \eqref{tildestrong}; i.e.
		$$\forall\;d\in\NN_{>0}\;\exists\;C\ge 1\;\forall\;j\in\NN:\;\;\;\widetilde{M}^d_j\le CN_j.$$
		
		\item[$(v)$] The weights satisfy $v_N\hyperlink{trianglecpow}{\vartriangleleft^{\mathfrak{c}}}v_M$.
	\end{itemize}
\end{proposition}

\demo{Proof}
$(i)\Rightarrow(ii)$ and $(v)\Rightarrow(i)$ are clear. $(iv)\Leftrightarrow(v)$ holds by $(c)$ in Lemma \ref{charactabstractweightcor1}, $(iv)\Rightarrow(iii)$ is clear and $(iii)\Rightarrow(iv)$ holds as follows: Let $d\in\NN_{>0}$ be given and take $d_1\in\NN_{>0}$ such that $d\le cd_1$ with $c$ appearing in \eqref{charactprop3equ1foropow}. Recall that we get $\widetilde{M}^d\le\widetilde{M}^{cd_1}$, see \cite[$(3.26)$, Lemma 3.19]{weightedentirecharacterization}, and then \eqref{charactprop3equ1foropow} yields the conclusion.\vspace{6pt}

It remains to show $(ii)\Rightarrow(iii)$ and we follow \cite[Prop. 3.18 $(i)$]{weightedentirecharacterization}. For this we are using the auxiliary sequence
\begin{equation}\label{auxilarysequ}
	\underline{N}^c_j:=\sup_{t\ge 0}\frac{t^j}{\exp(c\omega_N(t))},\;\;\;c\in\NN_{>0};
\end{equation}
see \cite[$(3.28)$]{weightedentirecharacterization}. Crucially, recall that $\underline{N}^c\in\hyperlink{LCset}{\mathcal{LC}}$ for any $c\in\NN_{>0}$ and, moreover, we have
\begin{equation}\label{goodequivalenceclassic}
	\forall\;c\in\NN_{>0}\;\exists\;D>0\;\forall\;t\ge 0:\;\;\;\omega_{\underline{N}^c}(t)\le c\omega_N(t)\le 2\omega_{\underline{N}^c}(t)+D;
\end{equation}
see \cite[$(3.29)$]{weightedentirecharacterization}.

Let now $c$ be the parameter from \eqref{charactprop3equforopow}. By \cite[Thm. 3.4 $(i)$]{weightedentirecharacterization} we know that $v_{\underline{N}^c}$ is essential and by applying the comments from Section \ref{essentialsection} to this weight there exists an optimal $f_{v_{\underline{N}^c}}\in H^{\infty}_{v_{\underline{N}^c}}(\CC)$, see \eqref{approximable}, such that $\frac{1}{v_{\underline{N}^c}}\hyperlink{sim}{\sim} t\mapsto M(f_{v_{\underline{N}^c}},t)$. By the first estimate in \eqref{goodequivalenceclassic} and the inclusion \eqref{charactprop3equforopow} one has $H^{\infty}_{v_{\underline{N}^c}}(\CC)\subseteq H^{\infty}_{v^c_N}(\CC)\subseteq H^{\infty}_{v^{c'}_M}(\CC)$ for all $c'>0$. We consider the function $f^2_{v_{\underline{N}^c}}$ and by gathering everything we arrive at
\begin{align*}
	&\exists\;c\in\NN_{>0}\;\exists\;A,D\ge 1\;\forall\;c'>0\;\exists\;C\ge 1\;\forall\;t\ge 0:
	\\&
	\frac{1}{AD}\exp(c\omega_N(t))\le\frac{1}{A}\exp(2\omega_{\underline{N}^c}(t))\le M(f^2_{v_{\underline{N}^c}},t)\le C\exp(2c'\omega_M(t)),
\end{align*}
for the first estimate the second part of \eqref{goodequivalenceclassic} has been used. Note that $M(f^2_{v_{\underline{N}^c}},t)=M(f_{v_{\underline{N}^c}},t)^2$.

Let now $d\in\NN_{>0}$ be given and set $c'=\frac{1}{2d}$. Then we use the previous estimate and \eqref{Prop32Komatsu} and obtain for all $j\in\NN$ that
\begin{align*}
	M_{cdj}&=M_{cj/(2c')}=\sup_{t\ge 0}\frac{t^{cj/(2c')}}{\exp(\omega_{M}(t))}=\left(\sup_{t\ge 0}\frac{t^{cj}}{\exp(2c'\omega_{M}(t))}\right)^{d}\le (ADC)^d\left(\sup_{t\ge 0}\frac{t^{cj}}{\exp(c\omega_{N}(t))}\right)^{d}
	\\&
	=(ADC)^d\left(\sup_{t\ge 0}\frac{t^j}{\exp(\omega_{N}(t))}\right)^{cd}=(ADC)^d(N_j)^{cd}.
\end{align*}
Thus \eqref{charactprop3equ1foropow} is verified.
\qed\enddemo

When combining Propositions \ref{charactprop2foropow} and \ref{charactprop3foropow} and Lemmas \ref{charactabstractweightcor} and \ref{charactabstractweightcor1} we get the following final statement:

\begin{theorem}\label{littleoweightsequcharactpow}
	Let $M,N\in\hyperlink{LCset}{\mathcal{LC}}$ be given. Then the following are equivalent:
	
	\begin{itemize}
		\item[$(i)$] The weight sequences satisfy \eqref{tildestrong}; i.e.
		$$\forall\;d\in\NN_{>0}\;\exists\;C\ge 1\;\forall\;j\in\NN:\;\;\;(M_{dj})^{1/d}=\widetilde{M}^d_j\le CN_j.$$
		
		\item[$(ii)$] The corresponding weights satisfy $v_N\hyperlink{trianglecpow}{\vartriangleleft^{\mathfrak{c}}}v_M$.
		
		\item[$(iii)$] \eqref{charactprop2equforopow} holds; i.e.
		$$\forall\;c,c'>0:\;\;\;H^{\infty}_{v^c_N}(\CC)\subseteq H^{\infty}_{v^{c'}_M}(\CC),$$
		with continuous inclusion.
		
		\item[$(iv)$] \eqref{charactprop3equforopow} holds; i.e.
		$$\exists\;c>0\;\forall\;c'>0:\;\;\;H^{\infty}_{v^c_N}(\CC)\subseteq H^{\infty}_{v^{c'}_M}(\CC),$$
		with continuous inclusion.
		
		\item[$(v)$] \eqref{charactprop2equ1foropow} holds; i.e.
		$$H^{\infty}_{\underline{\mathcal{N}}^{\mathfrak{c}}}(\CC)\subseteq H^{\infty}_{\overline{\mathcal{M}}^{\mathfrak{c}}}(\CC),$$
		with continuous inclusion.
	\end{itemize}
	If in addition either $M$ or $N$ has \hyperlink{mg}{$(\on{mg})$} and either $M$ or $N$ has \eqref{om1omegaMchar}, then the list of equivalences can be extended by
	\begin{itemize}
		\item[$(vi)$] The sequences satisfy $M\hyperlink{mtriangle}{\vartriangleleft}N$.
		
		\item[$(vii)$] The weights satisfy $v_N\hyperlink{trianglec}{\vartriangleleft_{\mathfrak{c}}}v_M$.
	\end{itemize}
\end{theorem}

\subsection{The dilatation-type weight function setting}

We start with the following technical statement.

\begin{lemma}\label{weightfctforolemma}
	Let $u$ and $w$ be normalized weight functions and consider the following assertions:
	\begin{itemize}
		\item[$(a)$] The associated sequences satisfy $M^{w}\hyperlink{mtriangle}{\vartriangleleft}M^{u}$.
		
		\item[$(b)$] The weights satisfy
		\begin{equation}\label{weightfctforolemmaequ}
			\exists\;B,C\ge 1\;\forall\;h>0\;\exists\;C_h\ge 1\;\forall\;t\ge 0:\;\;\;w(Bht)\le CC_h\sqrt{u(t)}.
		\end{equation}
		
		\item[$(c)$] The weights satisfy
		\begin{equation}\label{weightfctforolemmaequ1}
			\exists\;B_1,B_2\ge 1\;\forall\;h>0\;\exists\;C_h\ge 1\;\forall\;t\ge 0:\;\;\;w(B_1ht)\le B_2C_hu(t).
		\end{equation}
		
		\item[$(d)$] The weights satisfy $u\hyperlink{trianglec}{\vartriangleleft_{\mathfrak{c}}}w$.
	\end{itemize}
	
	Then we get:
	
	\begin{itemize}
		\item[$(*)$] If $u$ is convex, then $(a)\Rightarrow(b)$.
		
		\item[$(*)$] If either $u$ or $w$ has moderate growth, then $(b)\Rightarrow(c)$.
		
		\item[$(*)$] If either $u$ is convex and has moderate growth or $w$ is convex and has moderate growth, then $(c)\Rightarrow(d)$.
		
		\item[$(*)$] If $w$ is convex and either $u$ or $w$ has moderate growth, then $(d)\Rightarrow(a)$.
	\end{itemize}
\end{lemma}

Summarizing, if both $u$ and $w$ are convex and such that either $u$ or $w$ is of moderate growth, then all assertions listed before are equivalent.

\demo{Proof}
$(a)\Rightarrow(b)$ By definition we get
$$\forall\;h>0\;\exists\;C_h\ge 1\;\forall\;t\ge 0:\;\;\;\omega_{M^u}(t)\le\omega_{M^w}(ht)+\log(C_h)\Leftrightarrow v_{M^w,h}(t)\le C_hv_{M^u}(t).$$
Now use \eqref{omegavequiv}: We apply the second estimate there to $w$ and the first one to $u$ for which we require \emph{convexity.} Thus we arrive at
$$\exists\;A\ge 1\;\forall\;h>0\;\exists\;C_h\ge 1\;\forall\;t\ge 0:\;\;\;w(ht)\le v_{M^w}(ht)\le C_hv_{M^u}(t)\le AC_h\sqrt{u(t)},$$
which yields \eqref{weightfctforolemmaequ} with $B=1$ and $C=A$.\vspace{6pt}

$(b)\Rightarrow(c)$ If $u$ satisfies \eqref{om6forv}, then $\sqrt{u(Ht)}\le e^{H/2}u(t)$ for some $H\ge 1$ and all $t\ge 0$ and so
$$w(BHht)\le CC_h\sqrt{u(Ht)}\le CC_he^{H/2}u(t);$$
i.e. \eqref{weightfctforolemmaequ1} with $B_1=BH$ and $B_2=Ce^{H/2}$.

Similarly, if $w$ satisfies in addition \eqref{om6forv}, then by combining this with \eqref{weightfctforolemmaequ} we get
$$\exists\;B\ge 1\;\exists\;H\ge 1\;\forall\;h>0\;\exists\;C_h\ge 1\;\forall\;t\ge 0:\;\;\;w(HBht)\le e^{H}(w(Bht))^2\le e^H(CC_h)^2u(t);$$
i.e. \eqref{weightfctforolemmaequ1} with $B_1=BH$ and $B_2=C^2e^{H}$ and the choice $C_h^2$ for given $h$.

$(c)\Rightarrow(d)$ Let $c,d>0$ and assume that $u$ is convex and satisfies in addition \eqref{om6forv}. Then write $$\frac{w(t)}{u(ct/d)}=\frac{w(t)}{u((B_1h)^{-1}t)}\frac{u((B_1h)^{-1}t)}{u(ct/d)},$$
and apply \eqref{strongweightsrelationforv} to $u$: For given $d,c>0$ we have to choose $h>0$ small enough to ensure $(B_1h)^{-1}>Hc/d\Leftrightarrow\frac{d}{B_1Hc}>h$ with $H$ denoting the constant appearing in \eqref{om6forv} and $B_1$ the one appearing in \eqref{weightfctforolemmaequ1}. Thus $w(t)=o(u(ct/d))$ as $t\rightarrow+\infty$ follows and so we have shown $w_d(t)=o(u_c(t))$ for any $d,c>0$; i.e. $u\hyperlink{trianglec}{\vartriangleleft_{\mathfrak{c}}}w$. Note that in order to apply \eqref{strongweightsrelationforv} also convexity for $u$ is used.\vspace{6pt}

Similarly, if $u$ is arbitrary but $w$ satisfies in addition \eqref{om6forv} and is convex, then for all $h,d,c>0$ we write $$\frac{w_d(t)}{u_c(t)}=\frac{w(dt)}{u(ct)}=\frac{w(dt)}{w(B_1hct)}\frac{w(B_1hct)}{u(ct)}.$$
By \eqref{weightfctforolemmaequ1} and \eqref{strongweightsrelationforv} applied to $w$ we get $w(dt)=o(u(ct))$ for all $d>HB_1hc$. So, when given $d,c>0$, we are choosing again $h<\frac{d}{B_1Hc}$ in order to conclude.\vspace{6pt}

$(d)\Rightarrow(a)$ $u\hyperlink{trianglec}{\vartriangleleft_{\mathfrak{c}}}w$ means $w(ct)=o(u(t))$ as $t\rightarrow+\infty$ for all $c>0$. Then by the first estimate in \eqref{omegavequiv} applied to $w$, for which we require convexity, and the second one applied to $u$ we get $v^2_{M^w}(ct)=o(v_{M^u}(t))$ for all $c>0$; i.e. $v_{M^u}\hyperlink{trianglec}{\vartriangleleft_{\mathfrak{c}}}v^2_{M^w}$. Analogously as in $(b)\Rightarrow(c)$, if either $u$ or $w$ has property \eqref{om6forv}, then $v_{M^u}\hyperlink{trianglec}{\vartriangleleft_{\mathfrak{c}}}v_{M^w}$ follows; see also the proof of \cite[Prop. 3.22 $(ii)$]{weightedentirecharacterization}. Finally $(b)\Rightarrow(a)$ in Lemma \ref{charactabstractweightcor} finishes the conclusion.
\qed\enddemo

The following statement is analogous to Proposition \ref{charactprop2foro}.

\begin{proposition}\label{mixedweightstrong}
	Let $u$ and $w$ be weights. Consider the following assertions:
	\begin{itemize}
		\item[$(i)$] The weights satisfy $u\hyperlink{trianglec}{\vartriangleleft_{\mathfrak{c}}}w$.
		
		\item[$(ii)$] We have that
		\begin{equation}\label{mixedweightstrongequ}
			\forall\;c,c'>0:\;\;\;H^{\infty}_{u_c}(\CC)\subseteq H^{\infty}_{w_{c'}}(\CC),
		\end{equation}
		with continuous inclusion.
		
		\item[$(iii)$] We have that
		$$H^{\infty}_{\underline{\mathcal{U}}_{\mathfrak{c}}}(\CC)\subseteq H^{\infty}_{\overline{\mathcal{W}}_{\mathfrak{c}}}(\CC),$$
		with continuous inclusion.
	\end{itemize}
	Then $(i)\Rightarrow(ii)\Rightarrow(iii)$ holds.
\end{proposition}

\demo{Proof}
$(i)\Rightarrow(ii)$ Since $u\hyperlink{trianglec}{\vartriangleleft_{\mathfrak{c}}}w$ is equivalent to $u_d\hyperlink{mtriangle}{\vartriangleleft}w_c$ for all $c,d>0$, see the first part of \eqref{strongweightrelation}, and so the assertion follows immediately.

$(ii)\Rightarrow(iii)$ is clear.
\qed\enddemo

Next we deal with the converse statement which is analogous to Proposition \ref{charactprop3foro}; we make again a reduction to the weight sequence setting.

\begin{proposition}\label{mixedweightstrong1}
	Let $u$ and $w$ be normalized weights and consider the following inclusion relation (as sets):
	\begin{equation}\label{mixedweightstrong1equ0}
		\exists\;c>0\;\forall\;c'>0:\;\;\;H^{\infty}_{u,c}(\CC)\subseteq H^{\infty}_{w,c'}(\CC).
	\end{equation}
	\begin{itemize}
		\item[$(i)$] If $u$ is convex and \eqref{mixedweightstrong1equ0}, then
\begin{equation}\label{mixedweightstrong1equ}
		\exists\;A\ge 1\;\forall\;h>0\;\exists\;C_h\ge 1\;\forall\;t\ge 0:\;\;\;w(ht)\le AC_h\sqrt{u(t)},
\end{equation}
		which implies \eqref{weightfctforolemmaequ}.
		\item[$(ii)$] If either $u$ is convex and of moderate growth or if $w$ is convex and of moderate growth and \eqref{mixedweightstrong1equ} is valid, then
		\begin{equation}\label{mixedweightstrong1equ1}
			\exists\;B\ge 1\;\exists\;H\ge 1\;\forall\;h>0\;\exists\;C_h\ge 1\;\forall\;t\ge 0:\;\;\;w(Hht)\le BC_hu(t),
		\end{equation}
		with $H$ the constant appearing in \eqref{om6forv} for the particular weight.
		
		Finally, in this case \eqref{mixedweightstrong1equ1} implies $u\hyperlink{trianglec}{\vartriangleleft_{\mathfrak{c}}}w$.
	\end{itemize}
\end{proposition}

Note that \eqref{mixedweightstrongequ} is stronger than \eqref{mixedweightstrong1equ0}.

\demo{Proof}
$(i)$ Let $c>0$ be the fixed parameter appearing in \eqref{mixedweightstrong1equ0}. Then, by assumption and the second part of \eqref{omegavequiv} applied to $u$, we get
$$\exists\;c>0\;\forall\;c'>0:\;\;\;H^{\infty}_{v_{M^u,c}}(\CC)\subseteq H^{\infty}_{u,c}(\CC)\subseteq H^{\infty}_{w,c'}(\CC).$$
Hence, by applying this information to $\theta_{M^u,c}$ we have
$$\exists\;c>0\;\forall\;c'>0\;\exists\;C\ge 1\;\forall\;t\ge 0:\;\;\;\omega_{M^u}(ct/2)\le\log(C)-\log(w(c't))=\log(C)+\omega^{w}(c't),$$
which gives for all $j\in\NN$:
\begin{align*}
	M^{w}_j&=\sup_{t>0}\frac{t^j}{\exp(\omega^{w}(t))}=\sup_{t>0}\frac{(c't)^j}{\exp(\omega^{w}(c't))}\le C(c')^j\sup_{t>0}\frac{t^j}{\exp(\omega_{M^{u}}(ct/2))}
	\\&
	=C(2c'/c)^j\sup_{s>0}\frac{s^j}{\exp(\omega_{M^{u}}(s))}=C(2c'/c)^jM^u_j.
\end{align*}
Now let $c'\rightarrow 0$ and since $c>0$ is fixed we have verified $M^w\hyperlink{mtriangle}{\vartriangleleft}M^u$. Then $(a)\Rightarrow(b)$ in Lemma \ref{weightfctforolemma}, for which we require convexity for $u$, yields the conclusion. Note: When working with the optimal $f_{v_{M^u,c}}$ instead of $\theta_{M^u,c}$ then we can only get rid of the constant $2$ above but which does not effect the conclusion.\vspace{6pt}

$(ii)$ This follows from the proof of $(b)\Rightarrow(c)\Rightarrow(d)$ in Lemma \ref{weightfctforolemma}.
\qed\enddemo

Gathering Lemmas \ref{charactabstractweightcor}, \ref{weightfctforolemma} and Propositions \ref{mixedweightstrong} and \ref{mixedweightstrong1} we get the following characterization:

\begin{theorem}\label{littleoweightfctcharact}
	Let $u$ and $w$ be normalized weights. Consider the following assertions:
	\begin{itemize}
		\item[$(i)$] The associated weight sequences satisfy $M^{w}\hyperlink{mtriangle}{\vartriangleleft}M^{u}$.
		
		\item[$(ii)$] The weights satisfy $v_{M^u}\hyperlink{trianglec}{\vartriangleleft_{\mathfrak{c}}}v_{M^w}$.
		
		\item[$(iii)$] The weights satisfy $u\hyperlink{trianglec}{\vartriangleleft_{\mathfrak{c}}}w$.
		
		\item[$(iv)$] We have \eqref{mixedweightstrongequ}; i.e.
		$$\forall\;c,c'>0:\;\;\;H^{\infty}_{u_c}(\CC)\subseteq H^{\infty}_{w_{c'}}(\CC),$$
		with continuous inclusion.
		
		\item[$(v)$] We have \eqref{mixedweightstrong1equ0}; i.e.
		$$\exists\;c>0\;\forall\;c'>0:\;\;\;H^{\infty}_{u,c}(\CC)\subseteq H^{\infty}_{w,c'}(\CC),$$
		with continuous inclusion.
		
		\item[$(vi)$] We have
		$$H^{\infty}_{\underline{\mathcal{U}}_{\mathfrak{c}}}(\CC)\subseteq H^{\infty}_{\overline{\mathcal{W}}_{\mathfrak{c}}}(\CC),$$
		with continuous inclusion.
	\end{itemize}
	Then we get:
	\begin{itemize}
		\item[$(*)$] If $u$ is convex and of moderate growth, then $(i)\Leftrightarrow(ii)\Rightarrow(iii)\Leftrightarrow(iv)\Leftrightarrow(v)\Leftrightarrow(vi)$.
		
		\item[$(*)$] If both $u$ and $w$ are convex and either $u$ or $w$ is of moderate growth, then all listed assertions are equivalent.
	\end{itemize}
\end{theorem}

\subsection{The exponential-type weight function setting}

\begin{proposition}\label{mixedweightpowstrong}
	Let $u$ and $w$ be normalized weights. Consider the following assertions:
	\begin{itemize}
		\item[$(i)$] The weights satisfy $u\hyperlink{trianglec}{\vartriangleleft^{\mathfrak{c}}}w$.
		
		\item[$(ii)$] We have that
		\begin{equation}\label{mixedweightpowstrongequ}
			\forall\;c,d>0:\;\;\;H^{\infty}_{u^c}(\CC)\subseteq H^{\infty}_{w^d}(\CC),
		\end{equation}
		with continuous inclusion.
		\item[$(iii)$] We have that
		$$H^{\infty}_{\underline{\mathcal{U}}^{\mathfrak{c}}}(\CC)\subseteq H^{\infty}_{\overline{\mathcal{W}}^{\mathfrak{c}}}(\CC),$$
		with continuous inclusion.
	\end{itemize}
	Then $(i)\Rightarrow(ii)\Rightarrow(iii)$ holds.
\end{proposition}

\demo{Proof}
$(i)\Rightarrow(ii)$ Since $u\hyperlink{trianglec}{\vartriangleleft^{\mathfrak{c}}}w$ is equivalent to $u^d\hyperlink{mtriangle}{\vartriangleleft}w^c$ for all $c,d>0$, see the second part of \eqref{strongweightrelation}, and so the assertion follows immediately.

$(ii)\Rightarrow(iii)$ is clear.
\qed\enddemo

\begin{proposition}\label{mixedweightpowstrong1}
	Let $u$ and $w$ be normalized weights and assume that $u$ is convex. If (as sets)
	\begin{equation}\label{mixedweightpowstrong1equ}
		\exists\;c\in\NN_{>0}\;\forall\;d>0:\;\;\;H^{\infty}_{u^c}(\CC)\subseteq H^{\infty}_{w^d}(\CC),
	\end{equation}
	then $u\hyperlink{trianglecpow}{\vartriangleleft^{\mathfrak{c}}}w$ holds.
\end{proposition}

Note that \eqref{mixedweightpowstrong1equ} is weaker than \eqref{mixedweightpowstrongequ}.

\demo{Proof}
We follow the ideas given in the proof of $(ii)\Rightarrow(iii)$ in Proposition \ref{charactprop3foropow}; see also \cite[Prop. 3.27 $(i)$]{weightedentirecharacterization}:

Let $c$ be the parameter from \eqref{mixedweightpowstrong1equ}. We consider $\underline{M^u}^c$; i.e. the auxiliary sequence from \eqref{auxilarysequ} corresponding to the associated weight sequence $M^u$ from \eqref{vBMTweight1equ1}. Thus \eqref{approximable} holds for the corresponding class (recall Section \ref{essentialsection}) and so there exists $f_{v_{\underline{M^u}^c}}\in H^{\infty}_{v_{\underline{M^u}^c}}(\CC)$ such that $\frac{1}{v_{\underline{M^u}^c}}\hyperlink{sim}{\sim} t\mapsto M(f_{v_{\underline{M^u}^c}},t)$. The inclusion \eqref{mixedweightpowstrong1equ} together with \eqref{omegavequivnewnew} and \eqref{goodequivalenceclassic} give $H^{\infty}_{v_{\underline{M^u}^c}}(\CC)\subseteq H^{\infty}_{v^c_{M^u}}(\CC)\subseteq H^{\infty}_{u^c}(\CC)\subseteq H^{\infty}_{w^d}(\CC)$ for any $d>0$. Consider again the function $f^2_{v_{\underline{M^u}^c}}$ and so gathering all this information we arrive at
\begin{align*}
	&\exists\;c\in\NN_{>0}\;\exists\;A,B,D\ge 1\;\forall\;d>0\;\exists\;C\ge 1\;\forall\;t\ge 0:
	\\&
	\frac{1}{ADB^c}(u^{c/2}(t))^{-1}\le\frac{1}{AD}\exp(c\omega_{M^u}(t))\le\frac{1}{A}\exp(2\omega_{\underline{M^u}^c}(t))\le  M(f^2_{v_{\underline{M^u}^c}},t)\le C(w^d(t))^{-2}.
\end{align*}
The first estimate is valid by the first part of \eqref{omegavequiv} applied to $u$ and for which {\itshape convexity} is needed. The second estimate holds by the second part of \eqref{goodequivalenceclassic}. Summarizing, we obtain $(w^{4d/c}(t))\le(AB^cCD)^{2/c}u(t)$ and so, as $d\rightarrow 0$ and since $c$ is fixed, we have shown $u\hyperlink{ompreceq}{\preceq}w^{d_1}$ for all $d_1>0$ which is precisely \eqref{strongweightrelationpow}. By $(iii)$ in Remark \ref{strongconfusionrem} this fact is equivalent to $u\hyperlink{trianglecpow}{\vartriangleleft^{\mathfrak{c}}}w$.
\qed\enddemo

We summarize:

\begin{theorem}\label{littleoweightfctcharactpow}
	Let $u$ and $w$ be normalized weights and assume that $u$ is convex. Then the following are equivalent:
	\begin{itemize}
		\item[$(i)$] The weights satisfy $u\hyperlink{trianglec}{\vartriangleleft^{\mathfrak{c}}}w$.
		
		\item[$(ii)$] \eqref{mixedweightpowstrongequ} holds; i.e.
		$$\forall\;c,d>0:\;\;\;H^{\infty}_{u^c}(\CC)\subseteq H^{\infty}_{w^d}(\CC),$$
		with continuous inclusion.
		
		\item[$(iii)$] \eqref{mixedweightpowstrong1equ} holds; i.e.
		$$\exists\;c\in\NN_{>0}\;\forall\;d>0:\;\;\;H^{\infty}_{u^c}(\CC)\subseteq H^{\infty}_{w^d}(\CC),$$
		with continuous inclusion.
		
		\item[$(iv)$] We have
		$$H^{\infty}_{\underline{\mathcal{U}}^{\mathfrak{c}}}(\CC)\subseteq H^{\infty}_{\overline{\mathcal{W}}^{\mathfrak{c}}}(\CC),$$
		with continuous inclusion.
	\end{itemize}
	
	If $w$ is also convex, then the list of equivalences can be extended by:
	\begin{itemize}
		\item[$(v)$] The weights satisfy $v_{M^u}\hyperlink{trianglec}{\vartriangleleft^{\mathfrak{c}}}v_{M^w}$.
		
		\item[$(vi)$] The associated weight sequences satisfy \eqref{tildestrong}; i.e.
		$$\forall\;d\in\NN_{>0}\;\exists\;C\ge 1\;\forall\;j\in\NN:\;\;\;\widetilde{M^w}^d_j=(M^w_{dj})^{1/d}\le CM^u_j.$$
	\end{itemize}
\end{theorem}

\demo{Proof}
The equivalences of the first four assertions hold by Propositions \ref{mixedweightpowstrong} and \ref{mixedweightpowstrong1}.

$(i)\Leftrightarrow(v)$ is valid by the second part in \eqref{omegavequivnew}, i.e. $v_{M^u}\hyperlink{simpowc}{\sim^{\mathfrak{c}}}u$, $v_{M^w}\hyperlink{simpowc}{\sim^{\mathfrak{c}}}w$, and for this convexity is required for \emph{both} weights.

$(v)\Leftrightarrow(vi)$ is $(iv)\Leftrightarrow(v)$ in Proposition \ref{charactprop3foropow} applied to the associated weight sequences.
\qed\enddemo

\begin{remark}
\emph{The main statements from this section apply to the classes with $o$-growth restriction as well. More precisely, in view of Theorem \ref{oclassesnotnew} the exponential-type results Theorem \ref{littleoweightsequcharactpow} and Theorem \ref{littleoweightfctcharactpow} follow immediately and the same is valid for the dilatation-type weight sequence result Theorem \ref{littleoweightsequcharact}.}
	
\emph{Finally, the dilatation-type abstract weight function case requires a subtle change: In Proposition \ref{mixedweightstrong1} we have to consider $\theta_{M^u,d}$ for some (fixed) $0<d<c$ and with $c$ being the parameter from \eqref{mixedweightstrong1equ0}; for this recall \eqref{strongthetaequ}. Then follow the arguments there with $c$ being replaced by $d$ and so the conclusion follows.}
\end{remark}

\section{Comparison between dilatation-type and exponential-type systems}\label{compdilaexposection}
The goal of this section is to compare the inductive and projective structures defined by \eqref{weightsystems} and \eqref{weightsystems1} and to characterize equality as l.c.v.s. in terms of the given weights.

For this both technical growth restrictions \eqref{om6forv} and \eqref{om1forv}, equivalently expressed in terms of \hyperlink{om6}{$(\omega_6)$} and \hyperlink{om1}{$(\omega_1)$} for the function $\omega^v$ from \eqref{omegafromv}, are becoming crucial.

\subsection{On the sufficiency of conditions $(\omega_6)$ and $(\omega_1)$}
We start with the following:

\begin{proposition}\label{secondweighfctprop}
Let $v$ be a (normalized) weight function and let $\omega^v$ be the function from \eqref{omegafromv}.
\begin{itemize}
\item[$(i)$] If $v$ has \eqref{om1forv}, equivalently if $\omega^v$ has \hyperlink{om1}{$(\omega_1)$}, then
$$H^{\infty}_{\underline{\mathcal{V}}_{\mathfrak{c}}}(\CC)\subseteq H^{\infty}_{\underline{\mathcal{V}}^{\mathfrak{c}}}(\CC),\hspace{15pt}H^{\infty}_{\overline{\mathcal{V}}^{\mathfrak{c}}}(\CC)\subseteq H^{\infty}_{\overline{\mathcal{V}}_{\mathfrak{c}}}(\CC),$$
with continuous inclusions.

\item[$(ii)$] If $v$ has \eqref{om6forv}, equivalently if $\omega^v$ has \hyperlink{om6}{$(\omega_6)$}, then
$$H^{\infty}_{\underline{\mathcal{V}}^{\mathfrak{c}}}(\CC)\subseteq H^{\infty}_{\underline{\mathcal{V}}_{\mathfrak{c}}}(\CC),\hspace{15pt}H^{\infty}_{\overline{\mathcal{V}}_{\mathfrak{c}}}(\CC)\subseteq H^{\infty}_{\overline{\mathcal{V}}^{\mathfrak{c}}}(\CC),$$
with continuous inclusions.
\end{itemize}
Summarizing, if $v$ has both \eqref{om1forv} and \eqref{om6forv}, then as l.c.v.s.
$$H^{\infty}_{\underline{\mathcal{V}}_{\mathfrak{c}}}(\CC)=H^{\infty}_{\underline{\mathcal{V}}^{\mathfrak{c}}}(\CC),\hspace{15pt}H^{\infty}_{\overline{\mathcal{V}}_{\mathfrak{c}}}(\CC)=H^{\infty}_{\overline{\mathcal{V}}^{\mathfrak{c}}}(\CC).$$
\end{proposition}

\demo{Proof}
$(i)$ If $\omega^v$ has \hyperlink{om1}{$(\omega_1)$}, then (by iteration)
$$\forall\;C\ge 1\;\exists\;C'\ge 1\;\forall\;t\ge 0:\;\;\;\omega^v(Ct)\le C'\omega^v(t)+C'\Leftrightarrow v^{C'}(t)\le e^{C'}v_C(t),$$
yielding the first part. Similarly,
$$\forall\;0<c\le 1\;\exists\;0<c'\le 1\;\forall\;t\ge 0:\;\;\;\omega^v(t/c)\le\frac{1}{c'}\omega^v(t)+\frac{1}{c'}\Leftrightarrow v_c(t)\le ev^{c'}(t),$$
showing the second part.\vspace{6pt}

$(ii)$ If $\omega^v$ has \hyperlink{om6}{$(\omega_6)$}, then (by iteration)
$$\forall\;C\ge 1\;\exists\;C'\ge 1\;\forall\;t\ge 0:\;\;\;C\omega^v(t)\le\omega^v(C't)+C'\Leftrightarrow v_{C'}(t)\le e^{C'}v^C(t),$$
yielding the first part. Similarly,
$$\forall\;0<c\le 1\;\exists\;0<c'\le 1\;\forall\;t\ge 0:\;\;\;\frac{1}{c}\omega^v(t)\le\omega^v(t/c')+\frac{1}{c'}\Leftrightarrow v^c(t)\le e^{c/c'}v_{c'}(t),$$
showing the second part.
\qed\enddemo

When applying Proposition \ref{secondweighfctprop} to $v=v_M$, $\omega^v=\omega_M$ and by taking into account Lemma \ref{techlemma} we have:

\begin{corollary}\label{secondweighfctpropcor}
Let $M\in\hyperlink{LCset}{\mathcal{LC}}$ be given such that \hyperlink{mg}{$(\on{mg})$} and \eqref{om1omegaMchar} holds.

Then as l.c.v.s. we get
$$H^{\infty}_{\underline{\mathcal{M}}_{\mathfrak{c}}}(\CC)=H^{\infty}_{\underline{\mathcal{M}}^{\mathfrak{c}}}(\CC),\hspace{15pt}H^{\infty}_{\overline{\mathcal{M}}_{\mathfrak{c}}}(\CC)=H^{\infty}_{\overline{\mathcal{M}}^{\mathfrak{c}}}(\CC).$$
\end{corollary}
{\itshape Note:} By \cite[Prop. 3.4]{subaddlike} when assuming \hyperlink{mg}{$(\on{mg})$} then condition \eqref{om1omegaMchar} is equivalent to
$$\exists\;Q\in\NN_{>0}:\;\;\;\liminf_{j\rightarrow+\infty}\frac{\mu_{Qj}}{\mu_j}>1.$$

\subsection{On the necessity of conditions $(\omega_6)$ and $(\omega_1)$ for the weight sequence setting}\label{necsectionweightsequ}
The aim is now to prove the converse of Proposition \ref{secondweighfctprop} and Corollary \ref{secondweighfctpropcor}.

\begin{proposition}\label{secondweighfctpropconv}
Let $M\in\hyperlink{LCset}{\mathcal{LC}}$ be given.
\begin{itemize}
\item[$(i)$] The inclusion (as sets) $H^{\infty}_{\underline{\mathcal{M}}_{\mathfrak{c}}}(\CC)\subseteq H^{\infty}_{\underline{\mathcal{M}}^{\mathfrak{c}}}(\CC)$ implies \hyperlink{om1}{$(\omega_1)$} for $\omega_M$, whereas $H^{\infty}_{\underline{\mathcal{M}}^{\mathfrak{c}}}(\CC)\subseteq H^{\infty}_{\underline{\mathcal{M}}_{\mathfrak{c}}}(\CC)$ yields \hyperlink{om6}{$(\omega_6)$}.

\item[$(ii)$] The (continuous) inclusion $H^{\infty}_{\overline{\mathcal{M}}^{\mathfrak{c}}}(\CC)\subseteq H^{\infty}_{\overline{\mathcal{M}}_{\mathfrak{c}}}(\CC)$ implies \hyperlink{om1}{$(\omega_1)$} for $\omega_M$, whereas the (continuous) inclusion $H^{\infty}_{\overline{\mathcal{M}}_{\mathfrak{c}}}(\CC)\subseteq H^{\infty}_{\overline{\mathcal{M}}^{\mathfrak{c}}}(\CC)$ yields \hyperlink{om6}{$(\omega_6)$}.
\end{itemize}
\end{proposition}

\demo{Proof}
$(i)$ Crucially, we use $\theta_{M,c}$ and $\theta_M^c$ from Lemma \ref{charholomfctlemma}.

$H^{\infty}_{\underline{\mathcal{M}}_{\mathfrak{c}}}(\CC)\subseteq H^{\infty}_{\underline{\mathcal{M}}^{\mathfrak{c}}}(\CC)$ and \eqref{charholomfctlemmaequ} yield
$$\forall\;c>0\;\exists\;c'>0\;\exists\;D\ge 1\;\forall\;t\ge 0:\;\;\;\exp(\omega_M(ct/2))\le\theta_{M,c}(t)=|\theta_{M,c}(t)|\le D\exp(c'\omega_M(t)).$$
When taking $c=4$ we immediately get \hyperlink{om1}{$(\omega_1)$} for $\omega_M$ with $L:=\max\{c',\log(D)\}$.

Similarly, the inclusion $H^{\infty}_{\underline{\mathcal{M}}^{\mathfrak{c}}}(\CC)\subseteq H^{\infty}_{\underline{\mathcal{M}}_{\mathfrak{c}}}(\CC)$ and \eqref{charholomfctlemmaequ1} imply
$$\forall\;c\in\NN_{>0}\;\exists\;c'>0\;\exists\;D\ge 1\;\forall\;t\ge 0:\;\;\;\exp(c\omega_M(t/2^{1/c}))\le\theta^c_M(t)=|\theta^c_M(t)|\le D\exp(\omega_M(c't)).$$
We take $c=2$ and get $2\omega_M(t)\le\omega_M(c'\sqrt{2}t)+\log(D)$ for all $t\ge 0$; i.e. \hyperlink{om6}{$(\omega_6)$} with $H:=\max\{c'\sqrt{2},\log(D)\}$. Note that formally we only require the inclusion as sets.\vspace{6pt}

$(ii)$ We follow the proof of \cite[Prop. 3.14]{weightedentirecharacterization}; see also the citations there. Assume that $H^{\infty}_{\overline{\mathcal{M}}^{\mathfrak{c}}}(\CC)\subseteq H^{\infty}_{\overline{\mathcal{M}}_{\mathfrak{c}}}(\CC)$. Both spaces are Fr\'{e}chet and by the continuity of this inclusion we get
$$\forall\;d>0\;\exists\;c>0\;\exists\;C\ge 1\;\forall\;f\in H^{\infty}_{\overline{\mathcal{M}}^{\mathfrak{c}}}(\CC):\;\;\; \|f\|_{v_{M,d}}\le C\|f\|_{v^c_M}.$$
We apply this estimate to the family of monomials $(f_k)_{k\in\NN}$, so
$$\forall\;d>0\;\exists\;c>0\;\exists\;C\ge 1\;\forall\;k\in\NN\;\forall\;z\in\CC:\;\;\;\frac{|z|^k}{\exp(\omega_M(d|z|))}\le C\frac{|z|^k}{\exp(c\omega_M(|z|))}.$$
Now fix $d<1$ and w.l.o.g. $c':=\frac{1}{c}\in\NN_{>0}$ since $c\mapsto\|\cdot\|_{v^c_M}$ is non-increasing. Thus, by using \eqref{Prop32Komatsu} we obtain for all $j\in\NN$ that
\begin{align*}
M_{c'j}&=\sup_{t\ge 0}\frac{t^{c'j}}{\exp(\omega_{M}(t))}=\left(\sup_{t\ge 0}\frac{t^{j}}{\exp(c\omega_{M}(t))}\right)^{c'}\ge\frac{1}{C^{c'}}\left(\sup_{t\ge 0}\frac{t^{j}}{\exp(\omega_{M}(dt))}\right)^{c'}
\\&
=\frac{1}{C^{c'}}\frac{1}{d^{c'j}}\left(\sup_{s\ge 0}\frac{s^{j}}{\exp(\omega_{M}(s))}\right)^{c'}=\frac{1}{C^{c'}}\frac{1}{d^{c'j}}(M_j)^{c'}.
\end{align*}
Consequently, we have shown
\begin{equation}\label{om1omegaMchar1}
\exists\;c'\in\NN_{>0}\;\exists\;C\ge 1\;\forall\;j\in\NN_{>0}:\;\;\;\frac{(M_{c'j})^{1/(c'j)}}{(M_j)^{1/j}}\ge\frac{1}{C^{1/j}}\frac{1}{d};
\end{equation}
i.e. \eqref{om1omegaMchar} with $c'=L$. Thus \hyperlink{om1}{$(\omega_1)$} for $\omega_M$ follows by \cite[Thm. 3.1]{subaddlike}.\vspace{6pt}

Finally, if $H^{\infty}_{\overline{\mathcal{M}}_{\mathfrak{c}}}(\CC)\subseteq H^{\infty}_{\overline{\mathcal{M}}^{\mathfrak{c}}}(\CC)$ with continuous inclusion, then
$$\forall\;d>0\;\exists\;c>0\;\exists\;C\ge 1\;\forall\;f\in H^{\infty}_{\overline{\mathcal{M}}_{\mathfrak{c}}}(\CC):\;\;\; \|f\|_{v^d_M}\le C\|f\|_{v_{M,c}},$$
and applied to the monomials we get
$$\forall\;d>0\;\exists\;c>0\;\exists\;C\ge 1\;\forall\;k\in\NN\;\forall\;z\in\CC:\;\;\;\frac{|z|^k}{\exp(d\omega_M(|z|))}\le C\frac{|z|^k}{\exp(\omega_M(c|z|))}.$$
Take $d<1$ and assume w.l.o.g. $\frac{1}{d}\in\NN_{>0}$. Then for all $j\in\NN$
\begin{align*}
M_j&=\sup_{t\ge 0}\frac{t^{j}}{\exp(\omega_{M}(t))}=c^j\sup_{t\ge 0}\frac{t^{j}}{\exp(\omega_{M}(ct))}\ge\frac{1}{C}c^j\sup_{t\ge 0}\frac{t^{j}}{\exp(d\omega_{M}(t))}
\\&
=\frac{1}{C}c^j\left(\sup_{s\ge 0}\frac{s^{j/d}}{\exp(\omega_{M}(s))}\right)^{d}=\frac{1}{C}c^j(M_{j/d})^d.
\end{align*}
In particular, when choosing $d:=\frac{1}{2}$ then we have shown
\begin{equation}\label{mgdiag}
\exists\;c>0\;\exists\;C\ge 1\;\forall\;j\in\NN:\;\;\;M_{2j}\le C^2c^{-2j}(M_j)^2.
\end{equation}
By \cite[Thm. 9.5.1]{dissertation} applied to the constant matrix $\{M\}$, alternatively see also \cite[Thm. 1]{matsumoto}, the estimate \eqref{mgdiag} is equivalent to \hyperlink{mg}{$(\on{mg})$} for $M$ and so \cite[Prop. 3.6]{Komatsu73} yields \hyperlink{om6}{$(\omega_6)$} for $\omega_M$.
\qed\enddemo

By taking into account Lemma \ref{techlemma} and by combining Corollary \ref{secondweighfctpropcor} and Proposition \ref{secondweighfctpropconv} we get the following characterization:

\begin{theorem}\label{powerdilacharact}
Let $M\in\hyperlink{LCset}{\mathcal{LC}}$ be given. Then the following are equivalent:
\begin{itemize}
\item[$(i)$] $M$ has \eqref{om1omegaMchar} and \hyperlink{mg}{$(\on{mg})$}.

\item[$(ii)$] $\omega_{M}$ has \hyperlink{om1}{$(\omega_1)$} and \hyperlink{om6}{$(\omega_6)$}.

\item[$(iii)$] $v_M$ satisfies \eqref{om1forv} and \eqref{om6forv}.

\item[$(iv)$] As l.c.v.s. we get
$$H^{\infty}_{\underline{\mathcal{M}}_{\mathfrak{c}}}(\CC)=H^{\infty}_{\underline{\mathcal{M}}^{\mathfrak{c}}}(\CC).$$

\item[$(v)$] As l.c.v.s. we get
    $$H^{\infty}_{\overline{\mathcal{M}}_{\mathfrak{c}}}(\CC)=H^{\infty}_{\overline{\mathcal{M}}^{\mathfrak{c}}}(\CC).$$
\end{itemize}
\end{theorem}

\subsection{On the necessity of conditions $(\omega_6)$ and $(\omega_1)$ for the weight function setting}\label{necsectionweightfct}

We study the comparison between the structures from \eqref{weightsystems} and \eqref{weightsystems1} in the abstract weight function setting.

\begin{proposition}\label{secondweighfctpropconv1}
Let $u$ be a normalized and convex weight.
\begin{itemize}
\item[$(i)$] The inclusion (as sets) $H^{\infty}_{\underline{\mathcal{U}}_{\mathfrak{c}}}(\CC)\subseteq H^{\infty}_{\underline{\mathcal{U}}^{\mathfrak{c}}}(\CC)$ implies \eqref{om1forv} for $u$, whereas $H^{\infty}_{\underline{\mathcal{U}}^{\mathfrak{c}}}(\CC)\subseteq H^{\infty}_{\underline{\mathcal{U}}_{\mathfrak{c}}}(\CC)$ implies \eqref{om6forv}.

\item[$(ii)$] The (continuous) inclusion $H^{\infty}_{\overline{\mathcal{U}}^{\mathfrak{c}}}(\CC)\subseteq H^{\infty}_{\overline{\mathcal{U}}_{\mathfrak{c}}}(\CC)$ implies \eqref{om1forv} for $u$, whereas the (continuous) inclusion $H^{\infty}_{\overline{\mathcal{U}}_{\mathfrak{c}}}(\CC)\subseteq H^{\infty}_{\overline{\mathcal{U}}^{\mathfrak{c}}}(\CC)$ implies \eqref{om6forv}.
\end{itemize}
\end{proposition}

{\itshape Note:} The proof is reduced to the weight sequence setting by involving the associated weight sequence $M^u$ and for this convexity of $u$ is indispensable.

\demo{Proof}
$(i)$ By assumption and \eqref{omegavequivnewnew} we get
$$\forall\;c>0:\;\;\;H^{\infty}_{v_{M^u},c}(\CC)\subseteq H^{\infty}_{\underline{\mathcal{U}}_{\mathfrak{c}}}(\CC)\subseteq H^{\infty}_{\underline{\mathcal{U}}^{\mathfrak{c}}}(\CC).$$
These inclusions, Lemma \ref{charholomfctlemma} applied to $M^u$ and the first part of \eqref{omegavequivnew}, for which \emph{convexity} of $u$ is required, yield
\begin{align*}
&\exists\;A\ge 1\;\forall\;c>0\;\exists\;d>0\;\exists\;D\ge 1\;\forall\;t\ge 0:
\\&
\frac{1}{A}u^{-1/2}(ct/2))\le\exp(\omega_{M^u}(ct/2))\le\theta_{M^u,c}(t)=|\theta_{M^u,c}(t)|\le D(u^d(t))^{-1}.
\end{align*}
When taking $c=4$ then we get $u^{2d}(t)\le(AD)^2u(2t)=(AD)^2u_2(t)$ and thus \eqref{om1forv} is verified with $L:=\max\{2d,2\log(AD)\}$.

Similarly, by assumption and \eqref{omegavequivnewnew} we get
$$\forall\;c>0:\;\;\;H^{\infty}_{v^c_{M^u}}(\CC)\subseteq H^{\infty}_{\underline{\mathcal{U}}^{\mathfrak{c}}}(\CC)\subseteq H^{\infty}_{\underline{\mathcal{U}}_{\mathfrak{c}}}(\CC).$$
These inclusions, Lemma \ref{charholomfctlemma} applied to $M^u$ and the first part of \eqref{omegavequivnew} yield
\begin{align*}
&\exists\;A\ge 1\;\forall\;c\in\NN_{>0}\;\exists\;d>0\;\exists\;D\ge 1\;\forall\;t\ge 0:
\\&
\frac{1}{A^c}u^{-c/2}(t/2^{1/c}))\le\exp(c\omega_{M^u}(t/2^{1/c}))\le\theta^c_{M^u}(t)=|\theta^c_{M^u}(t)|\le D(u(dt))^{-1}.
\end{align*}
We choose $c=4$ and get $u_{d2^{1/4}}(t)=u_d(t2^{1/4})\le A^4Du^2(t)$. So \eqref{om6forv} is verified with $H:=\max\{d2^{1/4},\log(A^4D)\}$.\vspace{6pt}

$(ii)$ Assume now that $H^{\infty}_{\overline{\mathcal{U}}^{\mathfrak{c}}}(\CC)\subseteq H^{\infty}_{\overline{\mathcal{U}}_{\mathfrak{c}}}(\CC)$. Both spaces are Fr\'{e}chet and, similarly as before, by the continuity of this inclusion we get
$$\forall\;d>0\;\exists\;c>0\;\exists\;C\ge 1\;\forall\;f\in H^{\infty}_{\overline{\mathcal{U}}^{\mathfrak{c}}}(\CC):\;\;\; \|f\|_{u_d}\le C\|f\|_{u^c}.$$
We apply this estimate again to the family of monomials $(f_k)_{k\in\NN}$ and get
$$\forall\;d>0\;\exists\;c>0\;\exists\;C\ge 1\;\forall\;k\in\NN\;\forall\;z\in\CC:\;\;\;\frac{|z|^k}{\exp(\omega^u(d|z|))}\le C\frac{|z|^k}{\exp(c\omega^u(|z|))}.$$

Now fix $d<1$ and w.l.o.g. $c':=\frac{1}{c}\in\NN_{>0}$ since $c\mapsto\|\cdot\|_{u^c}$ is non-increasing. Then, by using \eqref{vBMTweight1equ1} and following the estimation in $(ii)$ in Proposition \ref{secondweighfctpropconv} we obtain \eqref{om1omegaMchar1} for $M^u$; i.e. \eqref{om1omegaMchar} for this sequence. Thus \hyperlink{om1}{$(\omega_1)$} for $\omega_{M^u}$ follows by \cite[Thm. 3.1]{subaddlike} and so $(II)$ in Lemma \ref{techlemma} yields the conclusion.\vspace{6pt}

Finally, if $H^{\infty}_{\overline{\mathcal{U}}_{\mathfrak{c}}}(\CC)\subseteq H^{\infty}_{\overline{\mathcal{U}}^{\mathfrak{c}}}(\CC)$ with continuous inclusion, then
$$\forall\;d>0\;\exists\;c>0\;\exists\;C\ge 1\;\forall\;f\in H^{\infty}_{\overline{\mathcal{U}}_{\mathfrak{c}}}(\CC):\;\;\; \|f\|_{u^d}\le C\|f\|_{u_c},$$
and this applied to the monomials gives
$$\forall\;d>0\;\exists\;c>0\;\exists\;C\ge 1\;\forall\;k\in\NN\;\forall\;z\in\CC:\;\;\;\frac{|z|^k}{\exp(d\omega^u(|z|))}\le C\frac{|z|^k}{\exp(\omega^u(c|z|))}.$$
Take $d<1$ such that $\frac{1}{d}\in\NN_{>0}$ and then, by \eqref{vBMTweight1equ1} and following the estimation in $(ii)$ in Proposition \ref{secondweighfctpropconv} we obtain with $d:=\frac{1}{2}$ that $M^u$ satisfies \eqref{mgdiag}. Thus \hyperlink{mg}{$(\on{mg})$} for $M^u$ is verified and by $(I)$ in Lemma \ref{techlemma} we are done.
\qed\enddemo

Propositions \ref{secondweighfctprop} and \ref{secondweighfctpropconv1} and Lemma \ref{techlemma} yield the following characterization:

\begin{theorem}\label{powerdilacharactforweights}
Let $u$ be a normalized and convex weight. Then the following are equivalent:
\begin{itemize}
\item[$(i)$] $u$ satisfies \eqref{om1forv} and \eqref{om6forv}.

\item[$(ii)$]  $\omega^u$ from \eqref{omegafromv} has \hyperlink{om1}{$(\omega_1)$} and \hyperlink{om6}{$(\omega_6)$}.

\item[$(iii)$] $\omega_{M^u}$ has \hyperlink{om1}{$(\omega_1)$} and \hyperlink{om6}{$(\omega_6)$}.

\item[$(iv)$] $M^u$ satisfies \eqref{om1omegaMchar} and \hyperlink{mg}{$(\on{mg})$}.

\item[$(v)$] As l.c.v.s. we get
$$H^{\infty}_{\underline{\mathcal{U}}_{\mathfrak{c}}}(\CC)=H^{\infty}_{\underline{\mathcal{U}}^{\mathfrak{c}}}(\CC).$$

\item[$(vi)$] As l.c.v.s. we get $$H^{\infty}_{\overline{\mathcal{U}}_{\mathfrak{c}}}(\CC)=H^{\infty}_{\overline{\mathcal{U}}^{\mathfrak{c}}}(\CC).$$
\end{itemize}
\end{theorem}

\begin{remark}
\emph{The main results Theorems \ref{powerdilacharact} and \ref{powerdilacharactforweights} are also valid for classes with $o$-growth condition:}

\emph{First, Proposition \ref{secondweighfctprop} clearly transfers to the $o$-type classes. Second, by Theorem \ref{oclassesnotnew} the weight sequence setting in Section \ref{necsectionweightsequ} is immediate since the $o$- and $O$-types coincide. Finally, concerning the general weight function case in Section \ref{necsectionweightfct} we point out that in view of \eqref{strongthetaequ} and \eqref{strongthetaequ1} part $(i)$ in Proposition \ref{secondweighfctpropconv1} follows; for this note that the crucial inclusions are valid for any parameter $c>0$. And part $(ii)$ in Proposition \ref{secondweighfctpropconv1} holds by using the same proof.}
\end{remark}

\appendix

\section{Erratum on closedness under point-wise multiplication}\label{erratum}
When treating $(ii)$ in Proposition \ref{secondweighfctpropconv}, more precisely when checking \eqref{mgdiag}, the author has seen that the analogous proofs of \cite[Prop. 4.7 \& 4.12]{weightedentirecharacterization} contain a gap. And a similar problem appears in \cite[Prop. 4.7]{Borelmapalgebraity}.\vspace{6pt}

We comment now in more detail and give correct arguments:

\begin{itemize}
\item[$(I)$] The mapping $f\mapsto f^2$ is not linear and hence the argument concerning the continuity given in \cite[Prop. 4.7]{weightedentirecharacterization} fails. The derived estimate \cite[$(4.8)$]{weightedentirecharacterization}, i.e.
    \begin{equation}\label{48equ}
    \exists\;A,B\ge 1\;\forall\;j\in\NN:\;\;\;M_{2j}\le AB^jM_j,
    \end{equation}
     is stronger than \eqref{mgdiag}, and so as \hyperlink{mg}{$(\on{mg})$}, since by assumption in \cite[Prop. 4.7]{weightedentirecharacterization} we have that $M$ is a log-convex weight sequence, see Definition \ref{defweightsequ}, hence $(M_j)^{1/j}\rightarrow+\infty$ and so $M_j\ge 1$ for all $j$ sufficiently large.

     In fact \eqref{48equ} is too strong, it yields the same conclusion as \cite[$(4.9)$]{weightedentirecharacterization} in the proof of \cite[Cor. 4.14]{weightedentirecharacterization} for any $M$ being a log-convex weight sequence: When iterating \eqref{48equ} then we get $M_{2^n}\le A^nB^{2^n-1}M_1$ for all $n\in\NN$ because the case $n=1$ yields $M_2\le ABM_1$ and $n\mapsto n+1$ holds because $M_{2^{n+1}}=M_{22^n}\le AB^{2^n}M_{2^n}\le AB^{2^n}A^nB^{2^n-1}M_1=A^{n+1}B^{2^{n+1}-1}M_1$. The case $n=0$ is trivial.

Thus $(M_{2^n})^{1/2^n}\le AB(M_1)^{1/2^n}$ for all $n\in\NN$ and, since $j\mapsto(M_j)^{1/j}$ is non-decreasing by log-convexity and $M_0=1$, for $k\in\NN_{>0}$ with $2^n\le k<2^{n+1}$ we obtain $(M_k)^{1/k}\le(M_{2^{n+1}})^{1/2^{n+1}}\le AB(M_1)^{1/2^{n+1}}$. Summarizing, $\sup_{j\in\NN_{>0}}(M_j)^{1/j}<+\infty$ is verified but this contradicts the basic requirements on $M$.

\item[$(II)$] For the sake of completeness we mention that \eqref{mgdiag}, equivalently \hyperlink{mg}{$(\on{mg})$} provided that $M$ is a log-convex weight sequence, yields by iteration the following estimate:
    $$\exists\;A\ge 1\;\forall\;n\in\NN:\;\;\;M_{2^n}\le A^{n2^{n-1}}(M_1)^{2^n},$$
    and so only $(M_{2^n})^{1/2^n}\le A^{n/2}M_1$. Again the case $n=0$ is trivial, the case $n=1$ holds by \eqref{mgdiag} and then iteration gives $M_{2^{n+1}}=M_{22^n}\le A^{2^n}(M_{2^n})^2\le A^{2^n}A^{2n2^{n-1}}(M_1)^{22^n}=A^{2^n(n+1)}(M_1)^{2^{n+1}}$.

\item[$(III)$] The correct proof of \cite[Prop. 4.7]{weightedentirecharacterization} is now as follows: By assumption the bilinear multiplication operator $\mathfrak{m}: (f,g)\mapsto f\cdot g$ is continuously acting on $H^{\infty}_{\overline{\mathcal{M}}_{\mathfrak{c}}}(\CC)$; recall the system \eqref{weightsequweightsystem}. Then we consider $\mathfrak{m}_{\Delta}: (f,f)\mapsto f^2$ and note that on a l.c.v.s. $E$ equipped with a system of continuous seminorms $\mathbb{P}:=\{p_i: i\}$ a bilinear mapping $T: E\times E\rightarrow E$ is continuous if and only if for each $p_0\in\mathbb{P}$ there exist $C>0$ and $p_1,p_2\in\mathbb{P}$ such that $p_0(T(x_1,x_2)) \le C p_1(x_1) p_2(x_2)$ for all $x_1, x_2 \in E$. We apply this to $T=\mathfrak{m}_{\Delta}$ and so \cite[$(4.7)$]{weightedentirecharacterization} turns into
    \begin{equation}\label{convolutorlemma3equ}
\forall\;d>0\;\exists\;c>0\;\exists\;C\ge 1\;\forall\;f\in H^{\infty}_{\overline{\mathcal{M}}_{\mathfrak{c}}}(\CC):\;\;\;\|f^2\|_{v_{M,d}}\le C\|f\|^2_{v_{M,c}}.
\end{equation}
When applying \eqref{convolutorlemma3equ} to the monomials $f_k=z^k$ we get
$$\forall\;d>0\;\exists\;c>0\;\exists\;C\ge 1\;\forall\;k\in\NN\;\forall\;z\in\CC:\;\;\;\frac{|z|^{2k}}{\exp(\omega_M(d|z|))}\le C\frac{|z|^{2k}}{\exp(2\omega_M(c|z|))},$$
and \eqref{Prop32Komatsu} gives for all $j\in\NN$ that
\begin{align*}
M_{2j}&=\sup_{t\ge 0}\frac{t^{2j}}{\exp(\omega_{M}(t))}=\sup_{s\ge 0}\frac{(ds)^{2j}}{\exp(\omega_{M}(ds))}\le Cd^{2j}\sup_{s\ge 0}\frac{s^{2j}}{\exp(2\omega_{M}(cs))}
\\&
=C\left(\frac{d}{c}\right)^{2j}\sup_{u\ge 0}\frac{u^{2j}}{\exp(2\omega_{M}(u))}=C\left(\frac{d}{c}\right)^{2j}(M_j)^2.
\end{align*}
Consequently, instead of \cite[$(4.8)$]{weightedentirecharacterization}, we obtain
\begin{equation*}
\exists\;A,B\ge 1\;\forall\;j\in\NN:\;\;\;M_{2j}\le AB^j(M_j)^2.
\end{equation*}
Thus \eqref{mgdiag} is verified and the arguments from $(ii)$ in Proposition \ref{secondweighfctpropconv} yield that \hyperlink{mg}{$(\on{mg})$} for $M$ holds.

\item[$(IV)$] Moreover, analogous comments apply to the proof of \cite[Prop. 4.12]{weightedentirecharacterization} verifying \eqref{mgdiag} and hence \hyperlink{mg}{$(\on{mg})$} for the associated weight sequence $M^u$.

\item[$(V)$] The arguments in the second part of the proof of \cite[Cor. 4.14]{weightedentirecharacterization} change as follows: We use the estimates from $(III)$ for $d=c=1$ since we are dealing here with the space $H^{\infty}_{v_M}(\CC)$. Then \cite[$(4.9)$]{weightedentirecharacterization} turns into the weaker condition
    $$\exists\;C\ge 1\;\forall\;j\in\NN:\;\;\;M_{2j}\le C(M_j)^2.$$
    However, proceeding as in $(I)$ or $(II)$, by iteration this estimate gives (with $A=C$)
    $$\exists\;A\ge 1\;\forall\;n\in\NN:\;\;\;M_{2^n}\le A^{2^n-1}(M_1)^{2^n}.$$
    The case $n=0$ is trivial and $n=1$ holds by taking $j=1$. Then iteration gives $M_{2^{n+1}}=M_{22^n}\le C(M_{2^n})^2\le CC^{2(2^n-1)}(M_1)^{2^{n+1}}=C^{2^{n+1}-1}(M_1)^{2^{n+1}}$ as desired. In \cite[Cor. 4.14]{weightedentirecharacterization} we have assumed that $M\in\hyperlink{LCset}{\mathcal{LC}}$ and hence the previous estimate gives $\sup_{j\in\NN_{>0}}(M_j)^{1/j}<+\infty$ analogously as shown in $(I)$; but this is a contradiction to the basic assumptions on $M$.

And the final comment in the proof of \cite[Cor. 4.14]{weightedentirecharacterization}, dealing with $M^u$, follows now by applying the previous explanations to $M^u$.

\item[$(VI)$] The same problem appears in the Beurling part of the proof of \cite[Prop. 4.7]{Borelmapalgebraity}; there working with weighted formal power series spaces, for unexplained notation please see \cite{Borelmapalgebraity}. By analogous reasons we obtain: On the right-hand side of \cite[$(4.8)$]{Borelmapalgebraity} the square of the corresponding seminorm is missing; this gives then the crucial estimate
    $$\exp(\omega_{\widehat{M}^{(\lambda)}}(s^2))\le C\exp(2\omega_{M^{(\kappa)}}(s/h_1)).$$
And then, one obtains by involving again \eqref{Prop32Komatsu}; i.e. \cite[$(4.7)$]{Borelmapalgebraity}, that $(M_j^{(\kappa)})^2\le Ch_1^{-2j}\widehat{M}^{(\lambda)}_j$. Thus $M^{(\kappa)}_{2j}$ has to be replaced by $(M_j^{(\kappa)})^2$. However, this is sufficient to verify \cite[$(4.6)$]{Borelmapalgebraity} and the argument given in the proof in \cite[Prop. 4.7]{Borelmapalgebraity} dealing with the clear estimate $(M_j^{(\kappa)})^2\le M^{(\kappa)}_{2j}$, which follows by log-convexity and normalization in view of \eqref{expsuperadd}, is becoming superfluous. Indeed, involving the wrong arguments we have shown a too strong estimate but \eqref{expsuperadd} has somehow hidden the mistake when proving \cite[$(4.6)$]{Borelmapalgebraity}.

These comments should be compared with the arguments and estimates in the proof of \cite[Thm. 5.1]{Borelmapalgebraity}.

Finally, in view of these observations we remark that for given $M\in\hyperlink{LCset}{\mathcal{LC}}$ in \cite[Prop. 5.2 $(5.6)$]{Borelmapalgebraity} we can replace $(m_j)^{2C}$ by $(m_{2j})^C$ with $m_j:=M_j/j!$. Thus \cite[$(5.6)$]{Borelmapalgebraity} is equivalent to
\begin{equation}\label{56alternative}
\exists\;C\in\NN_{>0}\;\exists\;D,h\ge 1\;\forall\;j\in\NN:\;\;\;m_{2j}\le Dh^j(m_{Cj})^{1/C}.
\end{equation}
This can be seen by a careful inspection of the proof of \cite[Prop. 5.2]{Borelmapalgebraity} or directly as follows: Obviously, \eqref{56alternative} implies \cite[$(5.6)$]{Borelmapalgebraity} because $j!^2(m_j)^2=(M_j)^2\le M_{2j}=(2j)!m_{2j}$ and by using Stirling's formula when changing $D$ and the geometric factor $h$ accordingly. Conversely, when applying \cite[$(5.6)$]{Borelmapalgebraity} to all even $j=2k$ then we get $m_{2k}\le D_1h_1^k(m_{2Ck})^{1/(2C)}$ for some $D_1,h_1\ge 1$ and all $k\in\NN$ and so \eqref{56alternative} holds with the choices $D_1,h_1$ and $C':=2C$.
\end{itemize}






\bibliographystyle{plain}
\bibliography{Bibliography}

\end{document}